\documentclass[preprint,12pt,sort&compress]{elsarticle}

\usepackage[margin=2.5cm]{geometry}
\usepackage{setspace}
\usepackage{lmodern}
\usepackage{amsmath,amssymb}
\usepackage{graphicx}
\usepackage{subfigure}
\usepackage{epsf,epsfig}
\usepackage{bm,bbm}
\usepackage{diagbox}
\usepackage{algpseudocode}
\usepackage{algorithmicx, algorithm}
\usepackage{physics}
\usepackage{color}
\usepackage{pxfonts}
\usepackage[footnotesize,bf]{caption}
\usepackage{mathtools}
\usepackage{hhline}
\usepackage[T1]{fontenc}
\usepackage{placeins}
\usepackage{multirow,dcolumn}
\usepackage{mathrsfs}
\usepackage{verbatim}
\usepackage{comment}
\usepackage{epstopdf}
\usepackage{lineno}
\usepackage[colorlinks=true,
linkcolor=blue,
anchorcolor=blue,
citecolor=blue,
urlcolor=blue,]{hyperref}
\usepackage[capitalise]{cleveref}
\usepackage{ulem}

\newtheorem{remark}{Remark}[section]

\graphicspath{{figures/}}
\journal{arXiv}

\begin{document}
\begin{frontmatter}

\title{Flow-based Bayesian filtering for high-dimensional nonlinear stochastic dynamical systems}

\author[TJU]{Xintong Wang}
\author[TJU]{Xiaofei Guan}
\author[SHNU]{Ling Guo}
\author[SJTU]{Hao Wu\corref{cor}}
\ead{hwu81@sjtu.edu.cn}
\address[TJU]{School of Mathematical Sciences, Tongji University, Shanghai 200092, China}
\address[SHNU]{Department of Mathematics, Shanghai Normal University, Shanghai 200234, China}
\address[SJTU]{School of Mathematical Sciences, Institute of Natural Sciences, and MOE-LSC, Shanghai Jiao Tong University, Shanghai 200240, China}
\cortext[cor]{Corresponding Author}  

\begin{abstract}
Bayesian filtering for high-dimensional nonlinear stochastic dynamical systems is a fundamental yet challenging problem in many fields of science and engineering. Existing methods face significant obstacles: Gaussian-based filters struggle with non-Gaussian distributions, while sequential Monte Carlo methods are computationally intensive and prone to particle degeneracy in high dimensions. Although generative models in machine learning have made significant progress in modeling high-dimensional non-Gaussian distributions, their inefficiency in online updating limits their applicability to filtering problems.
To address these challenges, we propose a flow-based Bayesian filter (FBF) that integrates normalizing flows to construct a novel latent linear state-space model with Gaussian filtering distributions. This framework facilitates efficient density estimation and sampling using invertible transformations provided by normalizing flows, and it enables the construction of filters in a data-driven manner, without requiring prior knowledge of system dynamics or observation models. Numerical experiments demonstrate the superior accuracy and efficiency of FBF.
\end{abstract}


\end{frontmatter}



\section{Introduction}\label{section:introduction}
Filtering problem is a crucial research area in data assimilation \cite{LKS:2015}, with applications that span numerical weather prediction, oceanography, signal processing, finance, etc. \cite{SKAMAROCK:JCP:2008, TANG:CMAME:2021, CHEN:JCP:2022}. The primary objective of filtering is to estimate the current state of a dynamical system by integrating noisy observations with prior model information. However, several significant challenges arise when applying filtering methods to more general and complex scenarios, such as nonlinearity and the curse of dimensionality. Developing an efficient computation model for accurate and stable state estimation is of substantial significance and has attracted great attention from scientists and engineers.

Classical filtering methods encompass variational methods and Bayesian filtering. Variational methods \cite{PLOGMANN:CMAME:2024, BENACEUR:CMAME:2024}, such as 3D-VAR \cite{MACK:CMAME:2020} and 4D-VAR \cite{CHANDRAMOULI:JCP:2020}, aim to find a best point-estimate of the system's state by minimizing the misfit for predicting measurements while penalizing deviations from the prior forecast. Bayesian filtering \cite{Särkkä:2013}, on the other hand, constructs the posterior distribution of the model state of interest through Bayesian principles, thereby effectively quantifying the uncertainties arising from prior information and observations. However, some well-known Bayesian techniques, such as the Kalman filter (KF) \cite{Kalman:JBE:1960} and its extensions, including the extended Kalman filter \cite{GONZALEZ:CMAME:2017} and the ensemble Kalman filter \cite{E:OD:2003, HJ:2018, BJ:JCP:2021a, DONOGHUE:CMAME:2022}, rely on a Gaussian ansatz, rendering them less suitable for applications involving non-Gaussian filtering distributions \cite{LJ:SJUQ:2024}. Therefore, many researchers utilized the sequential Monte Carlo method, referred to as the particle filter (PF), to obtain a more accurate estimate of the target posterior distribution \cite{GSS:IPFRSP:1993a, CCF:IEERSN:1999, PITT:JASA:1999}. Nevertheless, in cases where the system involves complex nonlinearity, PF may require a large number of particles to represent the state accurately, leading to high computational costs. An alternative to these methods is the quadrature Kalman filter (QKF) and its variants \cite{AH:ITSP:2008, JX:ACC:2019, NAIK:EJC:2023, LW:SIVP:2025}, which approximate nonlinear systems using quadrature rules. The QKF can offer more accurate estimates for non-Gaussian distributions while generally being less computationally expensive than PF. However, the computational cost of QKF increases significantly as the state dimension grows, making it less suitable for real-time applications in large-scale systems. To overcome this issue, Bao et al. proposed a novel score-based nonlinear filtering methodology \cite{BZZ:JCP:2024, BZZ:CMAME:2024}, which has achieved superior tracking accuracy on extremely high-dimensional stochastic dynamical systems.

The commonality among these aforementioned filters lies in their dependence on precise knowledge and accurate modeling of the underlying dynamics. In contrast, model-free data assimilation techniques, often referred to as data-driven methods, have garnered attention in scenarios where explicit dynamical models are unavailable. Instead, these techniques leverage simulation or experimental data that capture the system's dynamics to implement state estimation \cite{BRAJARD:JCS:2020, CSW:IP:2023, GHC:ITSP:2024, ROSAFALCO:CMAME:2024, WANG:JCAM:2025}. For example, recurrent neural network (RNN) and its variants, such as long short-term memory network (LSTM) \cite{SEPP:NC:1997, CHENG:CMAME:2024} and gated recurrent unit (GRU) \cite{JCKY:arxiv:2014} have demonstrated remarkable performance in time series-related tasks. However, they often struggle to capture the uncertainty \cite{BPGZTN:ICML:2019}. To this end, there have been significant advances in using probabilistic surrogate modeling techniques to perform state estimation with uncertainty quantification (UQ) \cite{KSBS:ICLR:2017, BPGZTN:ICML:2019, Yuta:arxiv:2024}. For example, a Kalman filtering-aided encoder-decoder model, termed the recurrent Kalman network (RKN), was proposed in \cite{BPGZTN:ICML:2019}. Similarly, Yuta et al. \cite{Yuta:arxiv:2024} introduced a variational autoencoder (VAE)-based Bayesian filter that integrates a latent linear Gaussian state transition model, enabling efficient inference of target filtering distributions. Despite these advancements, the aforementioned approaches rely on approximating the target distribution using tractable distributions (e.g., Gaussian distributions), which limits their applicability to scenarios involving complex filtering distributions.

In the past decade, normalizing flow (NF) based models have demonstrated great potential for estimating complex non-Gaussian distributions, as shown in the works of \cite{HZRS:2ICCVPRC:2016, DBLP:ICLR:2017, KDPDP:NIPS:2018, NOKW:SCI:2019, GWZ:JCP:2022, GWW:MLST:2024, wu2024reaction}. As a generative model, the NF transforms a simple probability distribution into a more complex one through a series of invertible and differentiable mappings, enabling both efficient density estimation and sampling. Notably, there has been growing interest in state estimation for nonlinear stochastic dynamical systems using NF models. For instance, Bézenac et al. \cite{BRBBK:NIPS:2020} proposed a novel approach for modeling complex multivariate time series by augmenting linear Gaussian state-space models with NFs, namely normalizing Kalman filter (NKF). Chen et al.~\cite{CWL:2I2ICIFF:2021} employed conditional normalizing flows (CNFs) that incorporate information from the latest observations, constructing more efficient proposal distributions for differentiable particle filters (DPF) \cite{JRB:RSSX:2018}, which uses differentiable programming to learn system dynamics and observation models directly from data. Nonetheless, several limitations persist within the current research. For example, the state variable in NKF is just a latent variable for forecasting, rendering the application of NKF infeasible for estimation of actual physical states in dynamical systems. The CNF-DPF method in \cite{CWL:2I2ICIFF:2021} struggles to facilitate efficient sampling in high-dimensional situations by following the particle filtering workflow.

In this paper, we put effort into establishing a novel approach to address nonlinear filtering problems, namely flow-based Bayesian filter (FBF). Specifically, we leverage normalizing flows to transform the original state and measurement variables into a latent space governed by a linear state-space model (SSM). This latent SSM consists of two key components: (1) a linear transition equation that describes the state update conditioned on the measurement, and (2) a linear observation equation that represents the one-step-ahead prediction of the measurement. Significantly, this formulation guarantees the Gaussianity of the filtering distribution at every time step, enabling efficient online calculation of its mean and covariance matrix through a Kalman filter-like procedure. Consequently, the target non-Gaussian filtering distribution can be efficiently evaluated and sampled by leveraging the invertibility of NFs. Compared to existing filtering methods (such as PF \cite{GSS:IPFRSP:1993a} and RKN \cite{BPGZTN:ICML:2019}), the proposed FBF enables accurate evaluation of high-dimensional non-Gaussian filtering distributions while maintaining significant computational efficiency.


The remainder of this paper is organized as follows. \cref{section:settings-preliminaries} outlines the statement of the nonlinear filtering problem and provides overviews of classical Bayesian filtering techniques and flow-based models. In \cref{section:FBF}, the framework of our FBF is introduced, where the algorithms for training and filtering approaches with FBF are described in detail. In \cref{section:examples}, we demonstrate the performance of the proposed method through several numerical examples. Finally, concluding remarks are provided in \cref{section:conclusion}.

\vspace{-2pt}
\section{Problem settings and preliminaries}\label{section:settings-preliminaries}
\subsection{Background}\label{subsection:Bayes-filter}
We consider a discrete-time stochastic dynamical system described by the following state-space model (SSM) \cite{GSS:IPFRSP:1993a}:
\begin{align*}
\bm{x}_{k} & =f(\bm{x}_{k-1},\bm{e}_{k}),\\
\bm{y}_{k} & =h(\bm{x}_{k},\bm{v}_{k}),
\end{align*}
where $k\in\mathbb{N}^+$ represents discrete time, $\bm{x}_k\in\mathbb{R}^{m}$ is the state vector, and $\bm{y}_k\in\mathbb{R}^{n}$ is the corresponding measurement. The state $\bm x_k$ evolves according to the nonlinear transition function $f$, while the measurement $\bm y_k$ is generated by the nonlinear measurement function $h$. The process noise $\bm e_k$ and measurement noise $\bm v_k$ are both independent and identically distributed.
For an SSM, the filtering problem involves estimating the probability distribution of the current state $\bm x_k$ given all available observations up to the present time, $p(\bm{x}_k|\bm{y}_{1:k})$, where $\bm{y}_{1:k}\triangleq\{\bm{y}_1,\bm{y}_2,\cdots,\bm{y}_k\}$. This estimation process is central to many applications, enabling the reconstruction of latent system dynamics from noisy measurements.

Bayesian filtering provides a general framework for recursively estimating the filtering distribution as new measurements become available. Using Bayes' rule, the filtering distribution at time $k$ can be expressed as:
\begin{eqnarray*}
p(\bm{x}_{k}|\bm{y}_{1:k}) & \propto & p(\bm{x}_{k},\bm{y}_{k}|\bm{y}_{1:k-1})\\
 & = & \int p(\bm{x}_{k},\bm{y}_{k}|\bm{x}_{k-1})p(\bm{x}_{k-1}|\bm{y}_{1:k-1})\mathrm{d}\bm{x}_{k-1},
\end{eqnarray*}
where $p(\bm{x}_{k-1}|\bm{y}_{1:k-1})$ represents the filtering distribution at the previous time step. Thus, the filtering distribution can be updated by recursively computing the above equation with each new observation.

When $f$ and $h$ are linear functions and the noise terms $\bm e_k$ and $\bm v_k$ follow Gaussian distributions, it can be shown that the filtering distribution at each time step is also Gaussian. In such cases, Bayesian filtering can be efficiently performed using the Kalman filter \cite{Kalman:JBE:1960}, which provides analytical updates for the mean and covariance of the filtering distribution. Moreover, for nonlinear SSMs, if the filtering distributions can be reasonably approximated by Gaussian distributions, several extended versions of the Kalman filter have been developed to address nonlinearity. A comprehensive review of these methods can be found in \cite{khodarahmi2023review}.

When filtering distributions deviate significantly from Gaussianity, the particle filter (PF) \cite{LKS:2015, GSS:IPFRSP:1993a} is an effective method for state estimation. It captures complex non-Gaussian posterior distributions through empirical approximations built from a set of random samples, known as particles. However, in many applications, particularly those involving high-dimensional SSMs, the finite samples often fail to accurately track the evolution of the filtering distributions. This limitation makes PF highly susceptible to particle degeneracy, a phenomenon where most particles carry negligible weights, thereby significantly reducing its effectiveness \cite{LJ:SJUQ:2024}. Bao et al. \cite{BZZ:JCP:2024, BZZ:CMAME:2024} address this issue by incorporating the score-based diffusion model with the Bayesian filtering framework. However, the aforementioned methods rely on the explicit formulation of the SSM, which may limit their effectiveness in scenarios where model information is incomplete or unavailable.

\subsection{Problem statement and research motivation}
This study addresses a challenging yet practical scenario: \textit{filtering in nonlinear systems where the SSM, defined by the transition and observation functions $f$ and $h$, is unknown.} Instead of relying on explicit mathematical representations, we assume access only to simulation or experimental data that capture the system's dynamics. The dataset typically consists of one or more trajectories of the state variables and their corresponding measurements. Such settings are common in real-world applications, including complex physical systems and data-driven scientific models, where deriving analytical formulations is often infeasible \cite{RPSDQ:NC:2024}.

To address these challenges, we utilize a specialized neural network framework, normalizing flows, to construct priors for the state and observation processes. Building on this foundation, we develop a filtering method that achieves computational efficiency and, even in high-dimensional settings, accurately models non-Gaussian state distributions.

\subsection{Normalizing flows}\label{subsection:NFs}
Before presenting our approach, we provide a brief review of normalizing flows (NFs).
NFs \cite{PGNE:JMLR:2021} are generative models that learn a complex target distribution by constructing an invertible transformation from a simple distribution using a sequence of differentiable mappings. Suppose $\bm x\in\mathbb{R}^{D}$ is a D-dimensional vector that follows a complex distribution. The NFs model enables establishment of an bijective transformation $\mathcal{T}$ such that $\bm z = \mathcal{T}(\bm x)\in\mathbb R^D$ is distributed according to a simple distribution $p_{\bm z}(\bm z)$ (i.e. Gaussian distribution), and the probability density function (PDF) of \( \bm{x} \) can be computed using the change of variables formula:
\[
p_{\bm{x}}(\bm{x}) = p_{\bm{z}}(\bm{z}) \left|\mathrm{det}\left(\frac{\partial \mathcal{T}(\bm{x})}{\partial \bm{x}}\right)\right|.
\]
For simplicity, we will omit the subscript notation in the PDFs when there is no ambiguity.
Additionally, for two random variables $\bm{x}$ and $\bm{x}'$, and their corresponding bijective transformations $\bm{z} = \mathcal{T}(\bm{x})$ and $\bm{z}' = \mathcal{T}'(\bm{x}')$ defined by NFs, we have
\begin{eqnarray}
p\left(\bm{x}|\bm{x}'\right) & = & \frac{p\left(\bm{x},\bm{x}'\right)}{p\left(\bm{x}'\right)}\nonumber\\
 & = & \frac{p\left(\bm{z},\bm{z}'\right)\cdot\left|\mathrm{det}\left(\frac{\partial\mathcal{T}(\bm{x})}{\partial\bm{x}}\right)\right|\cdot\left|\mathrm{det}\left(\frac{\partial\mathcal{T}'(\bm{x}')}{\partial\bm{x}'}\right)\right|}{p\left(\bm{z}'\right)\cdot\left|\mathrm{det}\left(\frac{\partial\mathcal{T}'(\bm{x}')}{\partial\bm{x}'}\right)\right|}\nonumber\\
 & = & p\left(\bm{z}|\bm{z}'\right)\cdot\left|\mathrm{det}\left(\frac{\partial\mathcal{T}(\bm{x})}{\partial\bm{x}}\right)\right|.\label{eq:conditional-nf}
\end{eqnarray}

Numerous variants of NFs models have been investigated in the literature, presenting a range of architectures and capabilities. One can refer to \cite{KPB:ITPAMI:2021} for a detailed review. In this work, we employ the RealNVP architecture that includes multiple affine coupling blocks \cite{DBLP:ICLR:2017} to construct the invertible transformation. 

\section{Flow-based Bayesian filter}\label{section:FBF}
\subsection{Model structure}\label{sec:model-structure}
Here we design the flow-based Bayesian filter (FBF), a novel method to address the challenges of filtering in nonlinear systems, particularly when the SSM is unknown, and the filtering distribution deviates significantly from Gaussianity. The core idea is to leverage NFs to map the original SSM into a latent space, where the state and measurement are transformed as 
\begin{equation}\label{eq:nf-transformation}
\bm{\chi}_{k} = \mathcal{T}(\bm{x}_{k}), \quad \bm{\gamma}_{k} = \mathcal{V}(\bm{y}_{k}),
\end{equation}
respectively. This transformation embeds the original SSM into the latent space, enabling the analytical update of the filtering distribution \( p(\bm{\chi}_k \mid \bm{\gamma}_{1:k}) \).

In this work, we formulate the latent SSM as:
\begin{eqnarray}
\bm{\chi}_{k} & = & A(\bm{\gamma}_{k}) + B(\bm{\gamma}_{k}) \bm{\chi}_{k-1} +  \bm{\epsilon}_{k},\quad  \bm{\epsilon}_{k}\sim\mathcal{N}(0,{Q}_{\chi}(\bm\gamma_k)),\label{eq:latent-update} \\
\bm{\gamma}_{k} & = & C + D \bm{\chi}_{k-1} + \bm{\nu}_{k},\quad \bm{\nu}_{k}\sim\mathcal{N}(0,{Q}_{\gamma}),\label{eq:latent-prediction}
\end{eqnarray}
where \( A(\cdot): \mathbb{R}^{n} \rightarrow \mathbb{R}^{m} \), \( B(\cdot): \mathbb{R}^{n} \rightarrow \mathbb{R}^{m \times m} \), and \( Q_{\chi}(\cdot): \mathbb{R}^{n} \rightarrow \mathbb{R}^{m \times m} \) are parameterized by deep neural networks (DNNs). The terms \( C \in \mathbb{R}^{n} \), \( Q_{\gamma} \in \mathbb{R}^{n \times n} \), and \( D \in \mathbb{R}^{n \times m} \) are trainable parameters, and \( \bm{\epsilon}_{k}, \bm{\nu}_{k} \) are i.i.d. latent noise. For simplicity, both \( Q_{\chi} \) and \( Q_{\gamma} \) are assumed to be diagonal matrices (see \ref{appendix:hyperparams} for modeling details).
\begin{figure}[ht]
	\centering
	\includegraphics[scale=0.60]{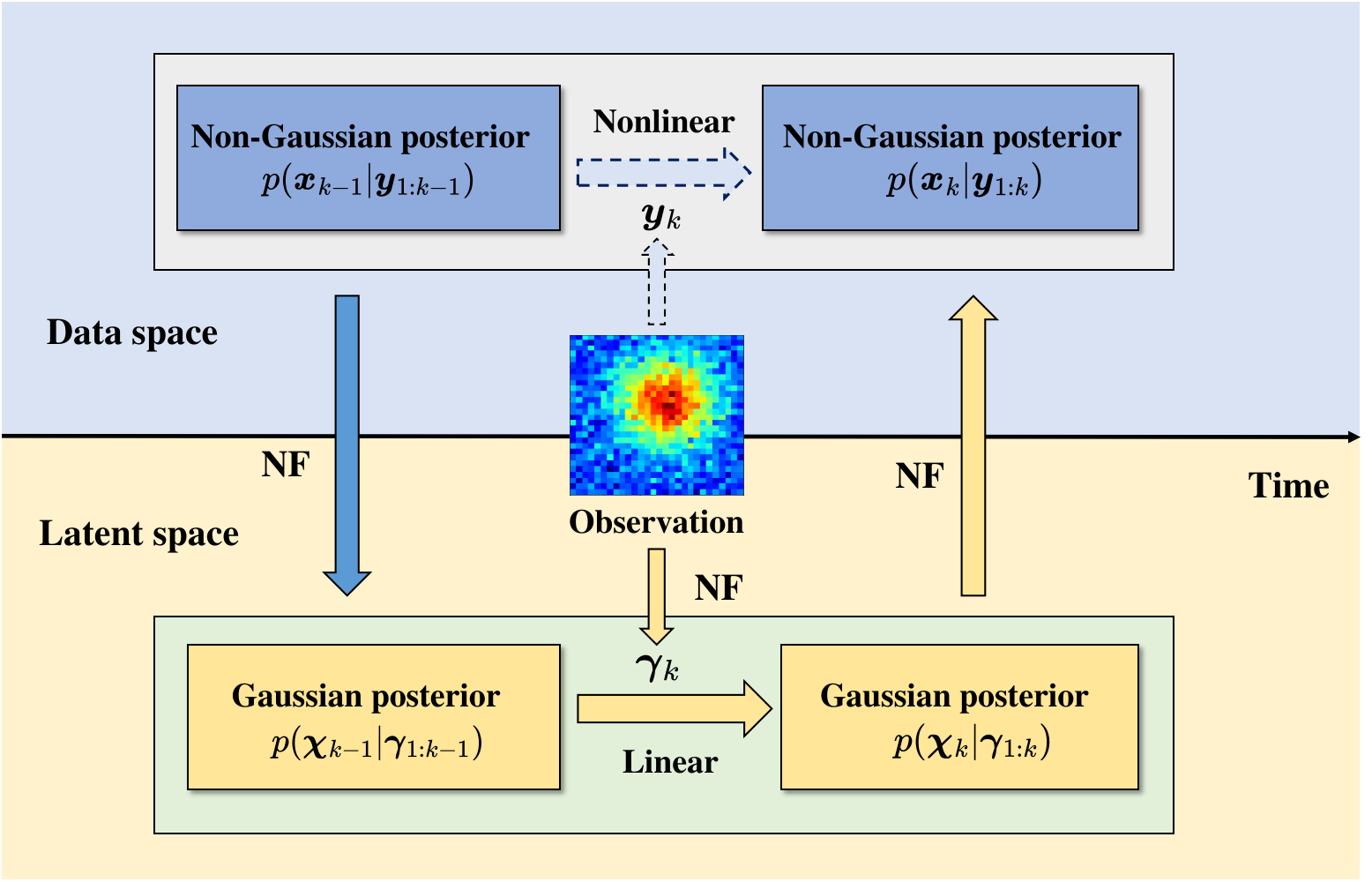}
	
	\vspace{+10pt}
	\caption{Schematic of FBF. The dashed arrows represent the execution of traditional Bayesian filtering, while the solid arrows indicate the direction of the proposed FBF.}
	\label{fig:FBF}
\end{figure}

We will prove in Section~\ref{sec:filtering-procedure} that the above model guarantees the filtering distribution \( p(\bm{\chi}_k \mid \bm{\gamma}_{1:k}) \) in the latent space remains Gaussian at every time step. Consequently, its mean and covariance matrix can be efficiently computed online using a procedure similar to the Kalman filter. Here, we provide an intuitive explanation. 
First, \eqref{eq:latent-prediction} represents the one-step-ahead prediction of the measurement. Since \(\bm{\gamma}_k\) is a linear combination of \(\bm{\chi}_{k-1}\) and Gaussian noise, with constant coefficient matrices \(C, D\), and the covariance \(Q_{\gamma}\) is independent of \(\bm{\chi}_{k-1}\), the posterior distribution of \(\bm{\chi}_{k-1}\) given \(\bm{\gamma}_k\) remains Gaussian if the filtering distribution \(p(\bm{\chi}_{k-1} \mid \bm{\gamma}_{1:k-1})\) at time \(k-1\) is Gaussian. 
Second, \eqref{eq:latent-update} describes the state update conditioned on the measurement. For a given measurement \(\bm{\gamma}_k\), the terms \(A, B, Q_{\chi}\) are constants. Therefore, \(\bm{\chi}_k\) can be considered as a linear function of two Gaussian random variables, \(\bm{\chi}_{k-1}\) and \(\bm{\epsilon}_k\), and the updated filtering distribution remains Gaussian.

While performing filtering in the latent space, the original-space filtering distribution can be reconstructed as follows:
\begin{equation}\label{eq:filtering-distribution-transformation}
p(\bm{x}_{k}\mid\bm{y}_{1:k})=p(\bm{\chi}_{k}\mid\bm{\gamma}_{1:k})\left|\mathrm{det}\left(\frac{\partial\mathcal{T}(\bm{x}_{k})}{\partial\bm{x}_{k}}\right)\right|.
\end{equation}

A systemic overview of the FBF is illustrated in \cref{fig:FBF}.

\begin{remark}
It is noteworthy that $A(\cdot)$, $B(\cdot)$, $C$ and $D$ are designed to be independent of $\bm\chi_{k-1}$, thereby preserving the linearity of the transformations in both \eqref{eq:latent-update} and \eqref{eq:latent-prediction}. Introducing dependency on $\bm\chi_{k-1}$ may render these transformations nonlinear, ultimately violating the Gaussianity of the latent SSM.
\end{remark}

\begin{remark}
The model defined by \eqref{eq:latent-update} and \eqref{eq:latent-prediction} differs slightly from the conventional linear SSM, particularly in the introduction of a nonlinear dependency on the observations in the state evolution equation \eqref{eq:latent-update}. It can be proved that conventional linear SSMs are a special case of our model, and our numerical experiments also demonstrate that the filter constructed using our model outperforms that based on standard linear SSM (see \ref{appendix:FBF-cSSM}).
\end{remark}



\subsection{Training Algorithm}\label{sec:training-process}
From \eqref{eq:conditional-nf}, the iterative equations for the latent SSM in \eqref{eq:latent-update} and \eqref{eq:latent-prediction} provide parametric models for the conditional probabilities \( p(\bm{x}_k \mid \bm{x}_{k-1}, \bm{y}_k) \) and \( p(\bm{y}_k \mid \bm{x}_{k-1}) \), as follows:
\begin{eqnarray}
p(\bm{x}_{k}|\bm{x}_{k-1},\bm{y}_{k}) & = & p(\bm{\chi}_{k}|\bm{\chi}_{k-1},\bm{\gamma}_{k})\left|\mathrm{det}\left(\frac{\partial\mathcal{T}(\bm{x}_{k})}{\partial\bm{x}_{k}}\right)\right|\nonumber \\
 & = & \mathcal{N}\left(\bm{\chi}_{k}|A(\bm{\gamma}_{k})+B(\bm{\gamma}_{k})\bm{\chi}_{k-1},Q_{\chi}(\bm{\gamma}_{k})\right)\left|\mathrm{det}\left(\frac{\partial\mathcal{T}(\bm{x}_{k})}{\partial\bm{x}_{k}}\right)\right|\nonumber \\
 & = & \mathcal{N}\left(\mathcal{T}(\bm{x}_{k})|A(\mathcal{V}(\bm{y}_{k}))+B(\mathcal{V}(\bm{y}_{k}))\mathcal{T}(\bm{x}_{k-1}),Q_{\chi}(\mathcal{V}(\bm{y}_{k}))\right)\left|\mathrm{det}\left(\frac{\partial\mathcal{T}(\bm{x}_{k})}{\partial\bm{x}_{k}}\right)\right|\nonumber \\
 & \triangleq & f_{s}(\bm{x}_{k-1},\bm{x}_{k},\bm{y}_{k}),\\
p(\bm{y}_{k}|\bm{x}_{k-1}) & = & p(\bm{\gamma}_{k}|\bm{\chi}_{k-1})\left|\mathrm{det}\left(\frac{\partial\mathcal{V}(\bm{y}_{k})}{\partial\bm{y}_{k}}\right)\right|\nonumber \\
 & = & \mathcal{N}\left(\bm{\gamma}_{k}|C+D\bm{\chi}_{k-1},Q_{\gamma}\right)\left|\mathrm{det}\left(\frac{\partial\mathcal{V}(\bm{y}_{k})}{\partial\bm{y}_{k}}\right)\right|\nonumber \\
 & = & \mathcal{N}\left(\mathcal{V}(\bm{y}_{k})|C+D\mathcal{T}(\bm{x}_{k-1}),Q_{\gamma}\right)\left|\mathrm{det}\left(\frac{\partial\mathcal{V}(\bm{y}_{k})}{\partial\bm{y}_{k}}\right)\right|\nonumber \\
 & \triangleq & f_{o}(\bm{x}_{k-1},\bm{y}_{k}),
\end{eqnarray}
where \( \mathcal{N}(\cdot \mid \bm{\mu}, \Sigma) \) denotes the PDF of the Gaussian distribution with mean \( \bm{\mu} \) and covariance matrix \( \Sigma \).

Given a training dataset consisting of a state trajectory \( \bm{x}_{0:K}^{\mathrm{train}} \) and the corresponding observations \( \bm{y}_{1:K}^{\mathrm{train}} \), we apply maximum likelihood estimation to train the NFs \( \mathcal{T}, \mathcal{V} \) and the latent SSM by maximizing the following objective function:
\[
\mathcal{L}(W) = \frac{\alpha}{K} \sum_{k=1}^{K} \log f_s(\bm{x}_{k-1}^{\mathrm{train}}, \bm{x}_k^{\mathrm{train}}, \bm{y}_k^{\mathrm{train}}) + \frac{\beta}{K} \sum_{k=1}^{K} \log f_o(\bm{x}_{k-1}^{\mathrm{train}}, \bm{y}_k^{\mathrm{train}}),
\]
where \( W \) represents all trainable parameters, and \( \alpha, \beta > 0\) are the weighting factors. In this work, we utilize stochastic gradient descent to optimize this objective (see Algorithm \ref{algorithm:FBF-train}). Furthermore, if the training data consists of multiple trajectories, a similar training objective can be used, which we omit for brevity.

In the modeling process outlined above, we have not considered the matching of the initial distributions in both the original and latent spaces. Noting that \( p(\bm{x}_0) = p(\bm{\chi}_0) \left| \frac{\partial \mathcal{T}(\bm{x}_0)}{\partial \bm{x}_0} \right| \), if \( p(\bm{x}_0) \) is given and its samples are available, we can add a term to the objective function:
\[
E_{p(\bm{x}_{0})}\left[\log\mathcal{N}(\mathcal{T}(\bm{x}_{0})|\bm{\mu}_{0},\bm{\Sigma}_{0})+\log\left|\mathrm{det}\left(\frac{\partial\mathcal{T}(\bm{x}_{0})}{\partial\bm{x}_{0}}\right)\right|\right],
\]
where \( \bm{\mu}_0 \) and \( \bm{\Sigma}_0 \) are also trainable parameters, ensuring that the initial state of the latent SSM follows a Gaussian distribution. However, in many applications, the initial state distribution \( p(\bm{x}_0) \) is simply chosen as a distribution with large variance, approximating an uninformative prior. In this work, the initial distribution is not considered during the FBF model training. Instead, after training, we estimate \( \bm{\mu}_0 \) and \( \bm{\Sigma}_0 \) by calculating the mean and covariance of the samples of \( \mathcal{T}(\bm{x}_0) \), and set \( p(\bm{\chi}_0) = \mathcal{N}(\bm{\chi}_0 | \bm{\mu}_0, \bm{\Sigma}_0) \). The proposed FBF, with this approximation, still achieves satisfactory filtering performance in our numerical examples.

\begin{algorithm}
	\caption{Training procedure for FBF}
	\vspace{+2pt}
	{\bf Input:}  Training data: $\{\bm{x}_{0:K}^{\mathrm{train}},\bm{y}_{1:K}^{\mathrm{train}}\}$, samples of the initial state: $\{\bm x_0^{(i)}\}_{i=1}^{N_0}$,
	number of epochs: $E_{\mathrm{train}}$, batch size: $N_{\mathrm{batch}}$, learning rate: $\eta$, weighting factors: $\alpha$ and $\beta$.
	\begin{algorithmic}
		\For{epoch = $1:E_{\mathrm{train}}$}
		\State Randomly draw indices $\{I_1,I_2,\ldots,I_{N_{\mathrm{batch}}}\}$ from the set $\{1,2,\ldots,K\}$;
            \State Calculate $$\hat{\mathcal{L}}=\frac{\alpha}{N_{\mathrm{batch}}}\sum_{i=1}^{N_{\mathrm{batch}}}\log f_{s}(\bm{x}_{I_{i}-1}^{\mathrm{train}},\bm{x}_{I_{i}}^{\mathrm{train}},\bm{y}_{I_{i}}^{\mathrm{train}})+\frac{\beta}{N_{\mathrm{batch}}}\sum_{k=1}^{N_{\mathrm{batch}}}\log f_{o}(\bm{x}_{I_{i}-1}^{\mathrm{train}},\bm{y}_{I_{i}}^{\mathrm{train}});$$
		\State Let $W\leftarrow W - \eta\nabla_{W}\hat{\mathcal{L}}$ to update all the involved parameters $W$;
		\EndFor
        \State Calculate $$
\bm{\mu}_{0}=\frac{\sum_{i=1}^{N_0}\mathcal{T}(\bm{x}_{0}^{(i)})}{N_0},\quad\bm{\Sigma}_{0}=\frac{\sum_{i=1}^{N_0}\left(\mathcal{T}(\bm{x}_{0}^{(i)})-\bm{\mu}_{0}\right)\left(\mathcal{T}(\bm{x}_{0}^{(i)})-\bm{\mu}_{0}\right)^{\top}}{N_0-1};
$$
	\end{algorithmic}
	\hspace*{0.001in} {\bf Output:}
	Trained FBF $\{\mathcal{T}, \mathcal{V}, A(\cdot), B(\cdot), Q_{\chi}(\cdot), C, D, Q_{\gamma},\bm\mu_0,\bm\Sigma_0\}$.
	\label{algorithm:FBF-train}
\end{algorithm}

\subsection{Filtering procedure}\label{sec:filtering-procedure}
We now explore how to solve nonlinear filtering problems using the FBF. Assume that at time step \( k-1 \), the latent filtering density is:
\[
p(\bm{\chi}_{k-1}|\bm{\gamma}_{1:k-1}) = \mathcal{N}(\bm{\chi}_{k-1}|\bm{\mu}_{k-1}, \bm{\Sigma}_{k-1}),
\]
and at time step \( k \), the new observation \( \bm{\gamma}_k = \mathcal{V}(\bm{y}_k) \) is available. The latent filtering distribution can be updated in the following two steps:

\paragraph{Step 1}
According to \eqref{eq:latent-prediction}, the posterior distribution of \( \bm{\chi}_{k-1} \) can be updated as:
\begin{equation}\label{equ:FBF_smooth}
\begin{aligned}
p(\bm{\chi}_{k-1}|\bm{\gamma}_{1:k}) & \propto\mathcal{N}(\bm{\chi}_{k-1}|\bm{\mu}_{k-1},\bm{\Sigma}_{k-1})\cdot\mathcal{N}(\bm{\gamma}_{k}|C+D\bm{\chi}_{k-1},Q_{\gamma})\\
 & \propto\exp\left(-\frac{(\bm{\chi}_{k-1}-\bm{\mu}_{k-1})^{\top}\bm{\Sigma}_{k-1}^{-1}(\bm{\chi}_{k-1}-\bm{\mu}_{k-1})}{2}\right)\\
 & \qquad\cdot\exp\left(-\dfrac{(\bm{\gamma}_{k}-C-D\bm{\chi}_{k-1})^{\top}Q_{\gamma}^{-1}(\bm{\gamma}_{k}-C-D\bm{\chi}_{k-1})}{2}\right)\\
 & \propto\mathcal{N}(\bm{\chi}_{k-1}|\bm{\mu}_{k-1}^{'},\bm{\Sigma}_{k-1}^{'}),
\end{aligned}
\end{equation}
where 
\begin{equation}\label{equ:FBF_smooth_mean_var}
\begin{aligned}
\bm{\mu}_{k-1}^{'} &= \bm{\mu}_{k-1} + \bm{\Sigma}_{k-1}D^{\top}Q_{\gamma}^{-1}(\bm{\gamma}_{k} - C - D\bm{\mu}_{k-1}),\\
\bm{\Sigma}_{k-1}^{'} &= \bm{\Sigma}_{k-1} - \bm{\Sigma}_{k-1}D^{\top}Q_{\gamma}^{-1}D\bm{\Sigma}_{k-1}.\\
\end{aligned}
\end{equation}

\paragraph{Step 2}
According to \eqref{eq:latent-update}, the latent filtering distribution at time step $k$ can be computed as:
\begin{equation}
\begin{aligned}
p(\bm{\chi}_{k}|\bm{\gamma}_{1:k}) &= \int p(\bm{\chi}_{k}|\bm{\gamma}_{k-1}, \bm{\gamma}_{k})~p(\bm{\chi}_{k-1}|\bm{\gamma}_{1:k})\mathrm{d}\bm{\chi}_{k-1}\\
&\propto \int \mathcal{N}(\bm{\chi}_{k}~|~A(\bm{\gamma}_{k}) + B(\bm{\gamma}_{k})\bm{\chi}_{k-1},Q_{\chi}(\bm{\gamma}_{k})~)\mathcal{N}(\bm{\chi}_{k-1}|\bm{\mu}_{k-1}^{'},\bm{\Sigma}_{k-1}^{'})\mathrm{d}\bm{\chi}_{k-1}\\
&= \mathcal{N}(\bm{\chi}_{k}|\bm{\mu}_{k},\bm{\Sigma}_{k}),
\end{aligned}
\end{equation}
where
\begin{equation}\label{equ:FBF_filter_mean_var}
\begin{aligned}
\bm{\mu}_{k} &= A(\bm{\gamma}_{k}) + B(\bm{\gamma}_{k})\bm{\mu}_{k-1}^{'}, \\
\bm{\Sigma}_{k} &= B(\bm{\gamma}_{k})\bm{\Sigma}_{k-1}^{'}B(\bm{\gamma}_{k})^{\top} + Q_{\chi}(\bm{\gamma}_{k}).
\end{aligned}
\end{equation}

Using \eqref{equ:FBF_smooth_mean_var} and \eqref{equ:FBF_filter_mean_var}, we can compute the latent filtering distribution efficiently online. Furthermore, by utilizing an invertible transformation \( \mathcal{T} \), we can perform efficient sampling of the target filtering distribution in the original space. Specifically, samples generated from \( p(\bm{\chi}_{k}|\bm{\gamma}_{1:k}) \) are transformed into the original space as follows:
\begin{equation}\label{equ:FBF_infer}
\bm{x}_k^{(j)} = \mathcal{T}^{-1}(\bm{\chi}_k^{(j)}), \quad \bm{\chi}_k^{(j)} \sim \mathcal{N}(\bm{\chi}_{k}|\bm{\mu}_{k},\bm{\Sigma}_{k}),
\end{equation}
where \( j \) denotes the sample index. The detailed procedure is outlined in \cref{algorithm:FBF-test}.

\begin{algorithm}[ht]
	\caption{Solving the nonlinear filtering problem by FBF}
	\hspace*{0.001in} {\bf Input:}
	Trained FBF: $\{\mathcal{T}, \mathcal{V}, A(\cdot), B(\cdot), Q_{\chi}(\cdot), C, D, Q_{\gamma}, \bm\mu_0, \bm\Sigma_0\}$, measurements $\{\bm{y}_1, \bm y_2, \ldots\}$, number of samples: $N_{\mathrm{samples}}$\\
	\vspace{-15pt}
	\begin{algorithmic}
		\For{$k$ = $1,2,\ldots$}
        \State Calculate $\bm{\mu}_{k}^{'}$ and $\bm{\Sigma}_{k}^{'}$ by \eqref{equ:FBF_smooth_mean_var};
        \State Calculate $\bm{\mu}_{k}$ and $\bm{\Sigma}_{k}$ by \eqref{equ:FBF_filter_mean_var};
		\For{$j$ = $1:N_{\mathrm{samples}}$}
		\State Draw $\bm{\chi}_{k}^{(j)}$ from $\mathcal{N}(\bm{\chi}_{k}|\bm{\mu}_{k},\bm{\Sigma}_{k})$;
		\State Calculate $\bm x_{k}^{(j)} = \mathcal{T}^{-1}(\bm{\chi}_{k}^{(j)})$;
		\EndFor
		\EndFor
	\end{algorithmic}
	\hspace*{0.001in} {\bf Output:} Samples $\{\bm x_k^{(j)}\}_{j=1}^{N_{\mathrm{samples}}}$ for each $k$.
	\label{algorithm:FBF-test}
\end{algorithm}

\section{Numerical results}\label{section:examples}
In this section, we present several numerical examples to demonstrate the effectiveness of the proposed method. First, we illustrate the performance of the approach using a two-dimensional nonlinear example, highlighting its ability to address non-Gaussian filtering problems. Next, we evaluate the performance of FBF in high-dimensional state estimation with the stochastic Lorenz-96 model. Finally, we apply FBF to a stochastic advection-diffusion model to validate its efficiency in solving filtering problems with nonlinear and sparse measurements. To benchmark the proposed method, we compare the performance of FBF with the following baselines:
\begin{itemize}
    \item Particle Filter (PF) \cite{GSS:IPFRSP:1993a}: Implements the standard sequential importance resampling.
    \item Recurrent Kalman Network (RKN) \cite{BPGZTN:ICML:2019}: Employs a recurrent neural network-based encoder-decoder architecture to obtain latent space representations of system states, with the latent state inferred using Kalman update equations. This method supports Gaussian or Bernoulli estimation of the filtering distribution in the original space\footnote{As all system states in this study are continuous variables, we use Gaussian-based decoders for the experiments.}.
    \item Conditional Normalizing Flow-based Differentiable Particle Filters (CNF-DPF) \cite{CWL:2I2ICIFF:2021}: Operates within the differentiable particle filter framework \cite{JRB:RSSX:2018}, employing NFs to directly learn the dynamics of nonlinear systems from data. Additionally, it leverages conditional NFs to construct proposal distributions for particle filtering.
\end{itemize}

Note that the PF utilizes true system SSMs during execution, while RKN and CNF-DPF, like FBF, rely on surrogate models derived from the same training data. The NFs used in our numerical examples are based on the RealNVP architecture \cite{DBLP:ICLR:2017}. Hyperparameter settings for both the network architecture and training for all methods are provided in \ref{appendix:hyperparams}. All computational simulations were performed on a server equipped with Intel Xeon E5-2650 v4 processors and NVIDIA GeForce RTX 2080 Ti GPUs.

To assess the performance of the proposed method, we use the following metrics (detailed computation is provided in \ref{appendix:metrics}):

\begin{itemize}
    \item Root Mean Squared Error (RMSE) \cite{HJ:2018}: Quantifies the accuracy of point estimates derived from the mean of the filtered distribution in recovering the true state values.
    \item Maximum Mean Discrepancy (MMD) \cite{Gretton:JMLR:2012}: Measures the accuracy of point estimates in a nonlinear feature space defined by a kernel function.
    \item Continuous Ranked Probability Score (CRPS) \cite{ding:arxiv:2024}: Assesses the consistency between the cumulative distribution functions of the inferred state variables and the true state values.
\end{itemize}

\subsection{Two-dimensional synthetic nonlinear model}\label{subsection:sinusoidal-example}
We consider a two-dimensional nonlinear SSM defined as follows \cite{RSNEvE:ITSP:2022}:
\begin{equation}\label{equ:ex1_equ_def}
\begin{split}
\bm{x}_{k} &= 0.9\,\mathrm{sin}(1.1\bm{x}_{k-1} + 0.1\pi) + 0.01 + \bm{e}_k,\quad \bm{e}_k\sim\mathcal{N}(0,q^2I_2),\ \  \bm{x}_k\in\mathbb{R}^2,\\
\bm{y}_{k} &= \mathrm{arctan}(\dfrac{x_{k,2}}{x_{k,1}}) +  \bm{v}_{k},\quad\quad \bm{v}_k\sim\mathcal{N}(0,r^2I_2),\ \  \bm{y}_{k}\in\mathbb{R}^2,
\end{split}
\end{equation}
where $\bm{x}_k = (x_{k,1}, x_{k,2})$ and $I_2\in\mathbb{R}^{2\times 2}$ denotes an identity matrix. We set a Gaussian prior for $\bm{x}_0$, where the initial mean is given by $(1, 1)\in\mathbb{R}^{2}$, and the covariance matrix is set to $0.1I_2\in\mathbb{R}^{2\times 2}$. In this experiment, we will verify the accuracy of our proposed methods under varying measurement noise levels. Thus we generate trajectory data of sequence length $K = 100$ with varying $r^2$: $r^2\in\{0.025, 0.05, 0.075, 0.1\}$. We set $q^2 = 0.1$, and use $N = 1000$ trajectories for training in each experiment.
\begin{figure}[!h]
	\centering
	\subfigure[]{
		\includegraphics[scale=0.43]{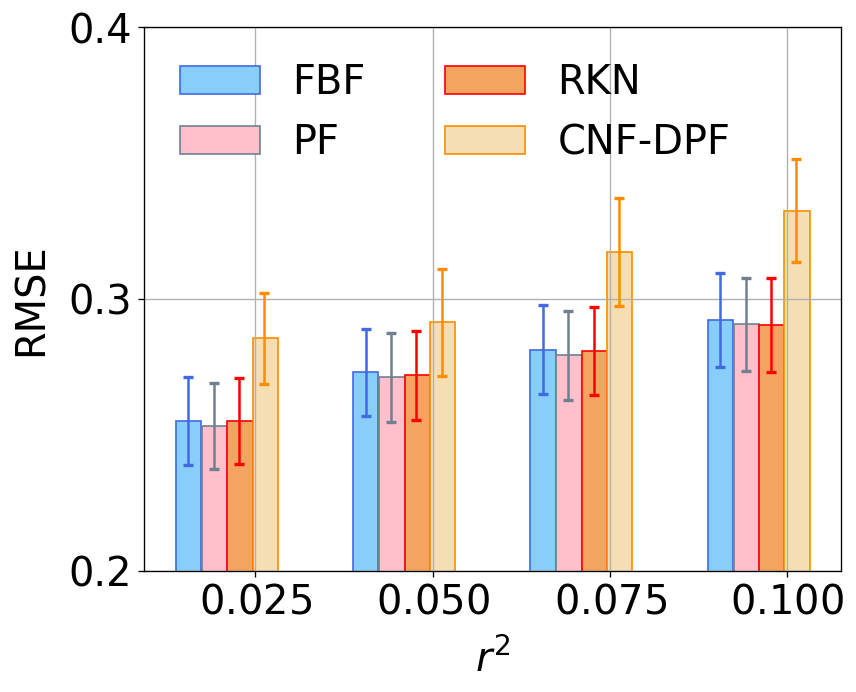}
	}
	\subfigure[]{
		\includegraphics[scale=0.43]{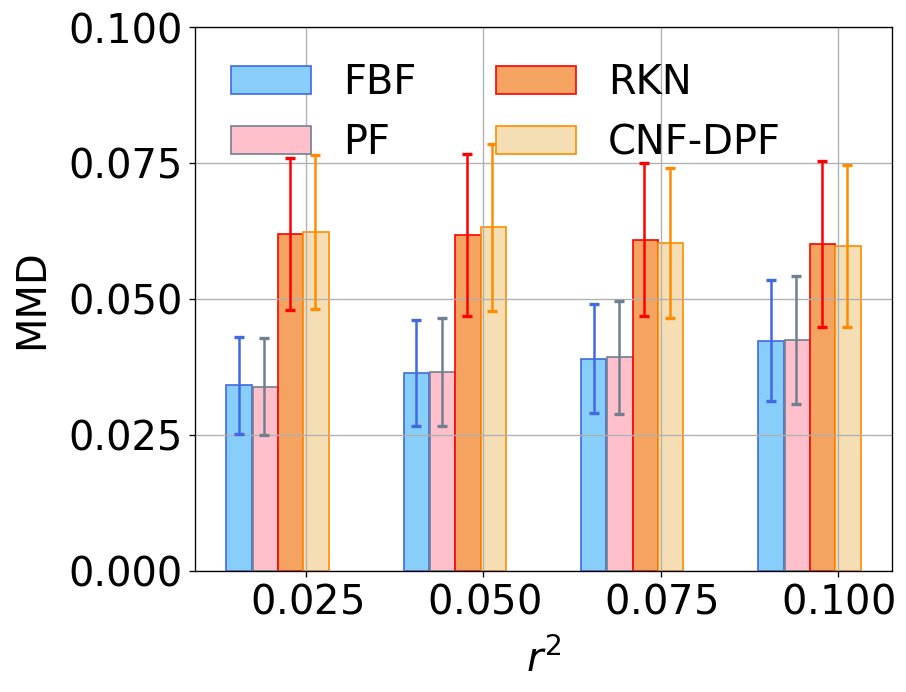}
	}\\
	\subfigure[]{
		\includegraphics[scale=0.43]{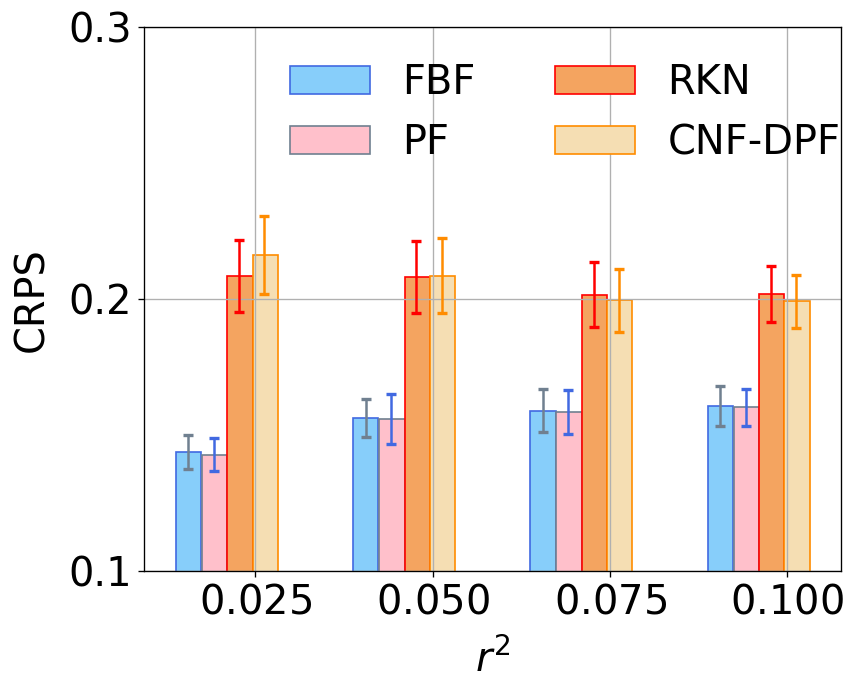}
	}
	\caption{Comparison of performance using different methods in \cref{subsection:sinusoidal-example}, where bars represent the mean and error bars indicate the standard deviation calculated over 200 test cases. (a) RMSE; (b) MMD; (c): CRPS.}
	\label{ex1_test_metric}
\end{figure}

We evaluate the performance of all methods using a test set comprising 200 trajectories and present the evaluation outcomes in \cref{ex1_test_metric}. As shown in \cref{ex1_test_metric}, for this two-dimensional example, the PF particles effectively track the evolution of the state distribution, achieving high filtering accuracy across various noise levels. In contrast, FBF, utilizing a surrogate model derived from the training data in the absence of knowledge of the system dynamics, yields results close to those of PF across all three metrics. Furthermore, when the noise level $r^2 \geq 0.05$, FBF significantly outperforms both RKN and CNF-DPF. These results indicate that the posterior distribution estimated by FBF closely matches the true distribution, demonstrating its effectiveness in probabilistic state estimation.

\begin{figure}[ht]
	\centering
	\subfigure{
		\includegraphics[scale=0.4]{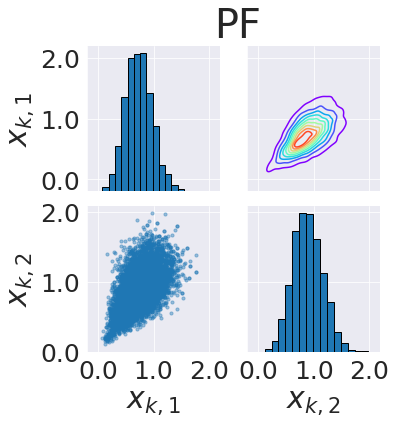}
	}
	\subfigure{
		\includegraphics[scale=0.4]{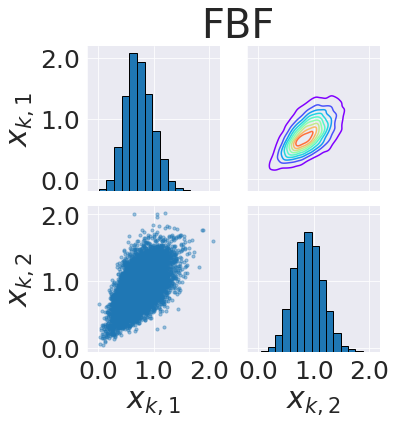}
	}\\
	\subfigure{
		\includegraphics[scale=0.4]{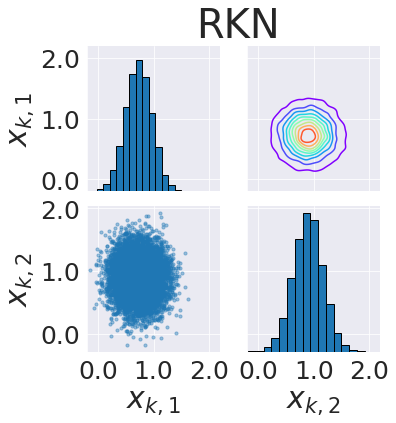}
	}
	\subfigure{
		\includegraphics[scale=0.4]{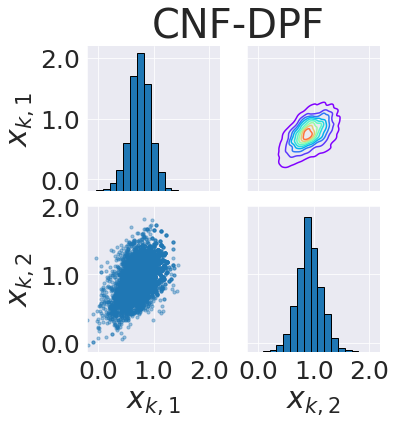}
	}
	\caption{The estimated filtering distribution of the inferred state for a given measurement at time step $k = 50$ by four methods. In each figure, the diagonal displays the histograms of state vector, while the lower left triangle and upper right triangle areas show the pairwise scatter diagram and density estimation of the same state vector, respectively.}
	\label{ex1_post_mesh}
\end{figure}

To visualize the effectiveness of the posterior approximation, we present the estimated filtering distributions for each method in \cref{ex1_post_mesh}. It is evident that the filtering distribution estimated by FBF closely matches that of the PF, which exhibits strong non-Gaussianity. The RKN method, due to its assumption of Gaussianity for the filtering distribution, fails to accurately capture the non-Gaussian state distribution. The filtering distribution provided by CNF-DPF is closer to that of PF, but some discrepancies remain.

The computational costs for each method, including both training and testing phases, are summarized in \cref{ex1_time_table}. It is evident that, with the same data size and number of training iterations, the training time required by FBF is significantly lower than that of the other two data-driven methods. During the testing phase, FBF requires more time than RKN due to the need to inversely transform samples from the latent space to the original data space. However, its computational cost remains considerably lower than that of CNF-DPF.


\begin{table}[!h]
	\centering
	\caption{Comparison of the computational cost required by different methods in Section \ref{subsection:sinusoidal-example}. For PF, the online time corresponds to the filtering process across $200$ test trajectories, whereas for the other data-driven methods, the offline time reflects $500$ training iterations for each method, and the online time represents the inference procedures for test cases, consistent with those of PF. The least time (in seconds) among the methods is printed in bold font.}
	\begin{tabular}{c|cccc}
		\hline
		Stage    & PF & FBF      & RKN        & CNF-DPF    \\ \hline
		Offline  & NA & $\bm{435.18}$   & $5965.41$ & $12840.28$ \\
		Online   & \bm{$17.09$} & $120.44$  & $52.49$    & $1084.45$  \\ \hline
	\end{tabular}
	\label{ex1_time_table}
\end{table}

\subsection{High-dimensional state estimation: Lorenz-96}\label{subsection:lorenz96_example}
Next, we demonstrate the performance of FBF for high-dimensional state estimation using the following stochastic Lorenz-96 model:
\begin{equation}\label{equ:lorenz96_example}
\begin{aligned}
\dfrac{\mathrm{d}x_{t,j}}{\mathrm{d}t} &= x_{t,j-1}(x_{t,j+1} - x_{t,j-2}) - x_{t,j} + F + \dfrac{\mathrm{d}W(t)}{\mathrm{d}t}, \quad j = 1,\cdots,m,\\
x_{0,j} &= \mathrm{sin}\big(\dfrac{2\pi j}{m}\big), \quad j = 1,\cdots,m,
\end{aligned}
\end{equation}
where $\bm{x}_t = (x_{t,1}, x_{t,2}, \cdots, x_{t,m})$ is an $m$-dimensional state vector, with boundary conditions $x_{t,-1} = x_{t,m-1}$, $x_{t,0} = x_{t,m}$, and $x_{t,m+1} = x_{t,1}$. The parameter $F = 8$ is a known forcing constant that induces chaotic behavior \cite{lorenz1998optimal}, and $W(t)$ represents an $m$-dimensional vector of independent standard Brownian motions.
The nonlinear measurement function from \cite{BZZ:JCP:2024} is adopted:
\begin{equation}
\bm y_t = \bm x_t^3 + \bm v_t,\quad\quad \bm v_t\sim\mathcal{N}(0,I_m).
\end{equation}
We simulate \eqref{equ:lorenz96_example} using a fourth-order Runge-Kutta scheme. During filtering, the observation interval is set to $\Delta t = 0.01$. In this case, we conduct experiments with varying dimensions: $m = 10,20,\cdots,50, 100$. For each $m$, we use $N = 2000$ labeled data points of length $K = 500$ for training and $N_{\mathrm{test}} = 200$ data points for testing. Note that numerical results for CNF-DPF are presented only for $m\leq 30$, as higher-dimensional scenarios require an excessive amount of GPU memory.

\begin{figure}[!h]
	\centering
	\subfigure[]{
		\includegraphics[scale=0.39]{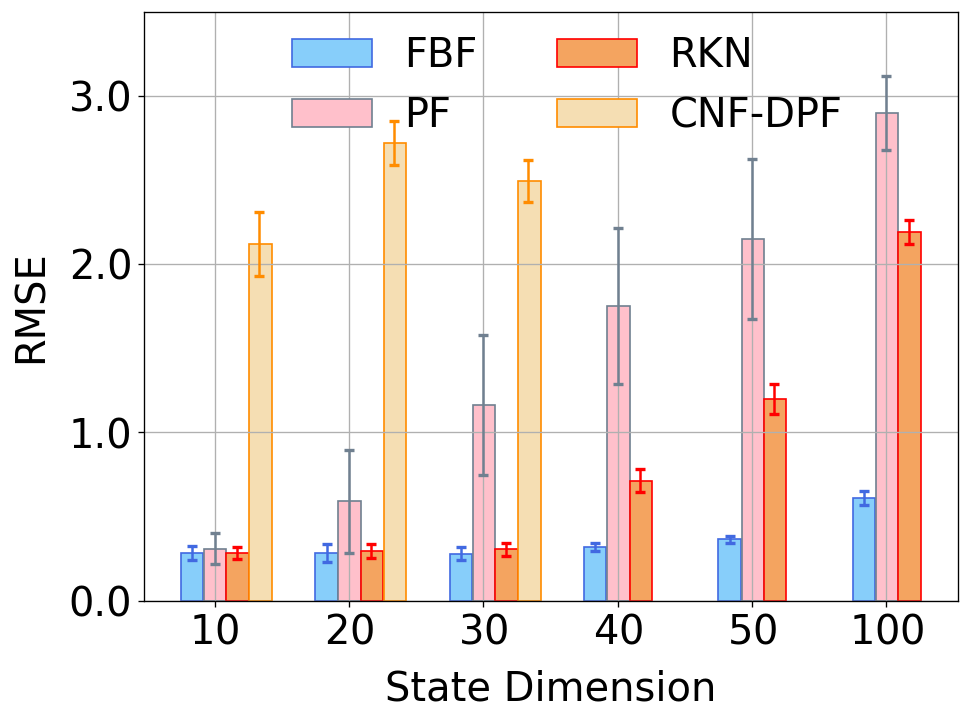}
	}
	\subfigure[]{
		\includegraphics[scale=0.39]{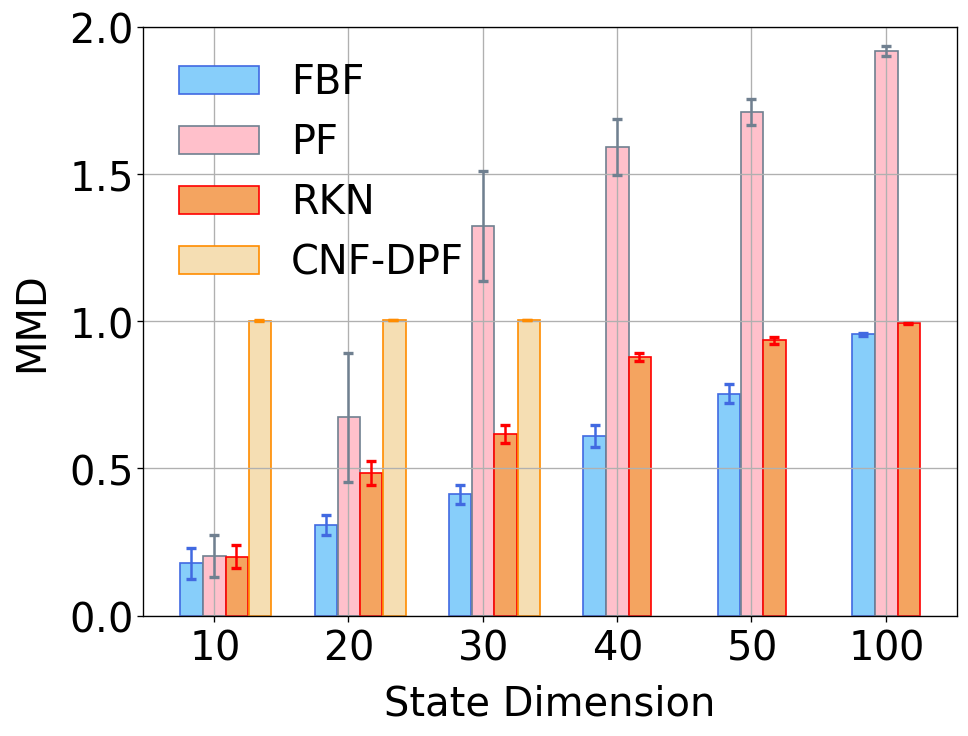}
	}
        \subfigure[]{
		\includegraphics[scale=0.39]{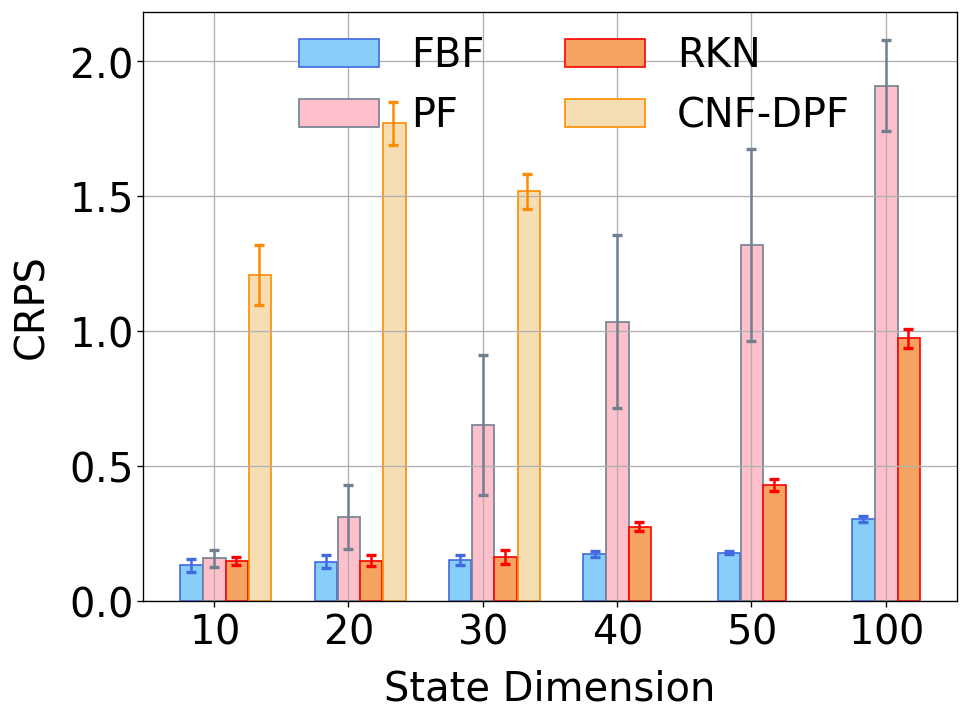}
	}
	\caption{Comparison of performance using different methods associated with varying dimensionality of states, where the bars represent the mean and the error bars indicate the standard deviation across 200 test cases. (a) RMSE; (b) MMD; (c) CRPS.}
	\label{fig:ex2_metrics}
\end{figure}

The performance comparison is presented in \cref{fig:ex2_metrics}. It is evident that the estimation errors of PF and CNF-DPF are significantly larger than those of the other two methods, indicating their limitations in high-dimensional scenarios. Additionally, while RKN performs comparably to FBF when $m \leq 30$, it struggles to accurately estimate states as dimensionality increases. In contrast, the proposed FBF demonstrates superior performance for $m \geq 40$. This discrepancy can be attributed to the simplicity of the RKN decoder, which is insufficient to capture complex state distributions in high-dimensional spaces. FBF, on the other hand, employs a more general latent state representation and a flexible NF-based architecture, allowing for efficient approximations of high-dimensional filtering distributions.

\begin{figure}[!h]
	\centering
	\subfigure{
		\includegraphics[scale=0.42]{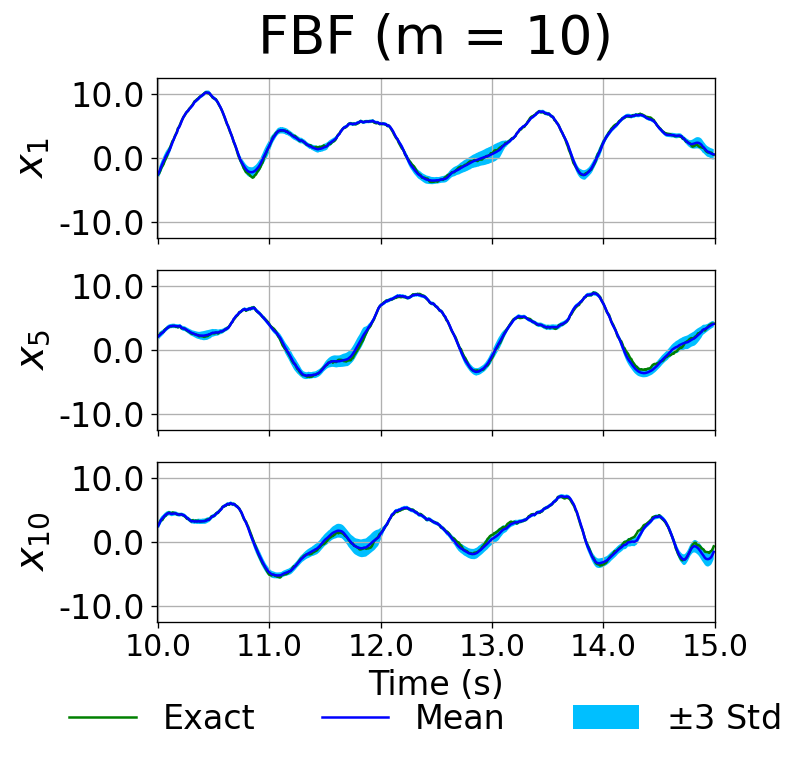}
	}
        \subfigure{
		\includegraphics[scale=0.42]{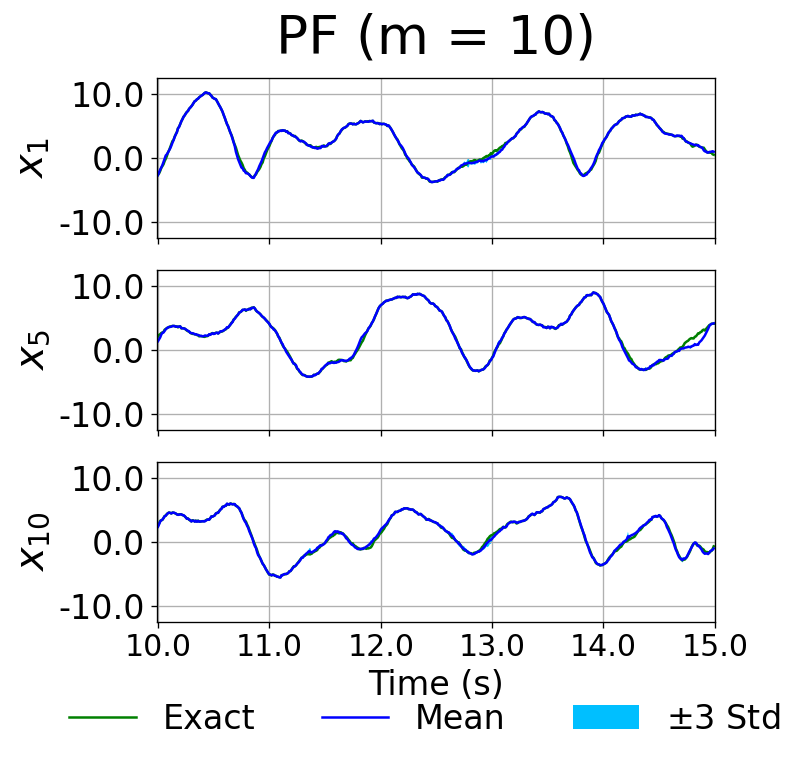}
	}
        \subfigure{
		\includegraphics[scale=0.42]{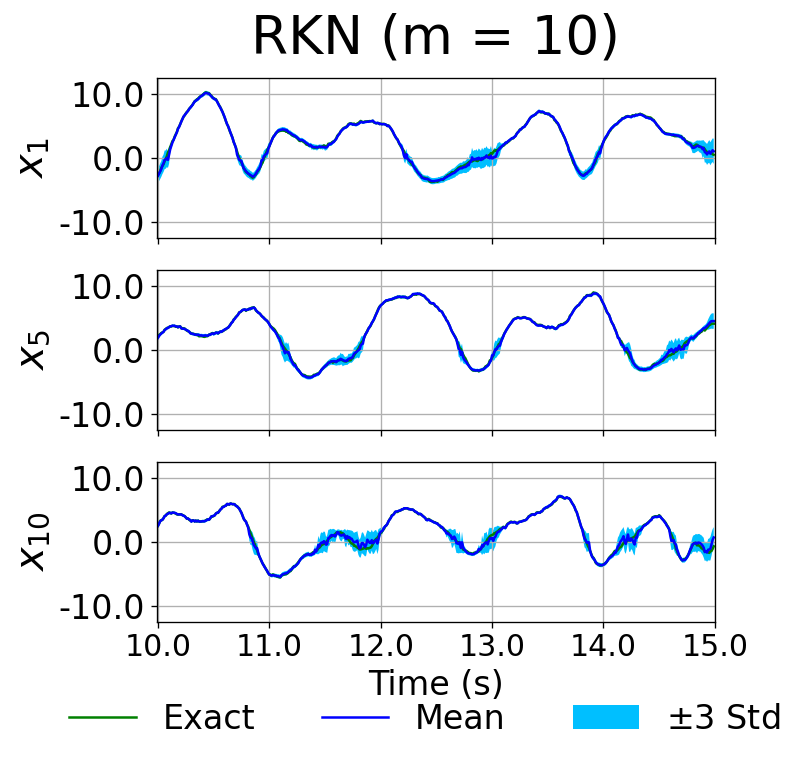}
	}
        \subfigure{
		\includegraphics[scale=0.42]{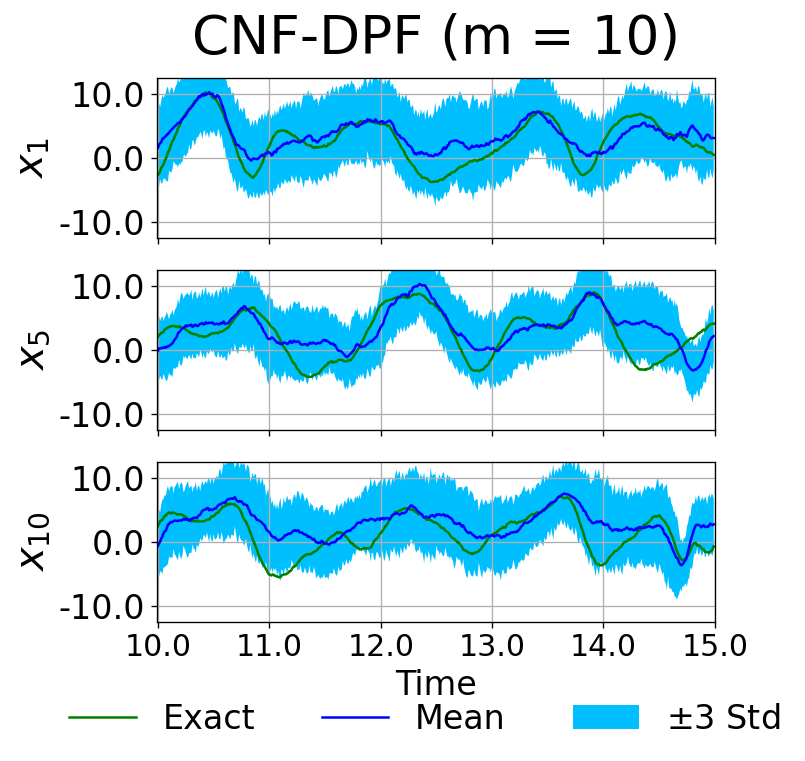}
	}
	\caption{Visualization of mean and uncertainty of the estimated posterior for a test measurement of $10$-dimensional Lorenz-96 systems. For clarity, we present the final 500 steps of the test trajectory of length $K_{\mathrm{test}} = 1500$.}
	\label{fig:ex2_post_interval_m=10}
\end{figure}

\begin{figure}[!h]
	\centering
	\subfigure{
		\includegraphics[scale=0.42]{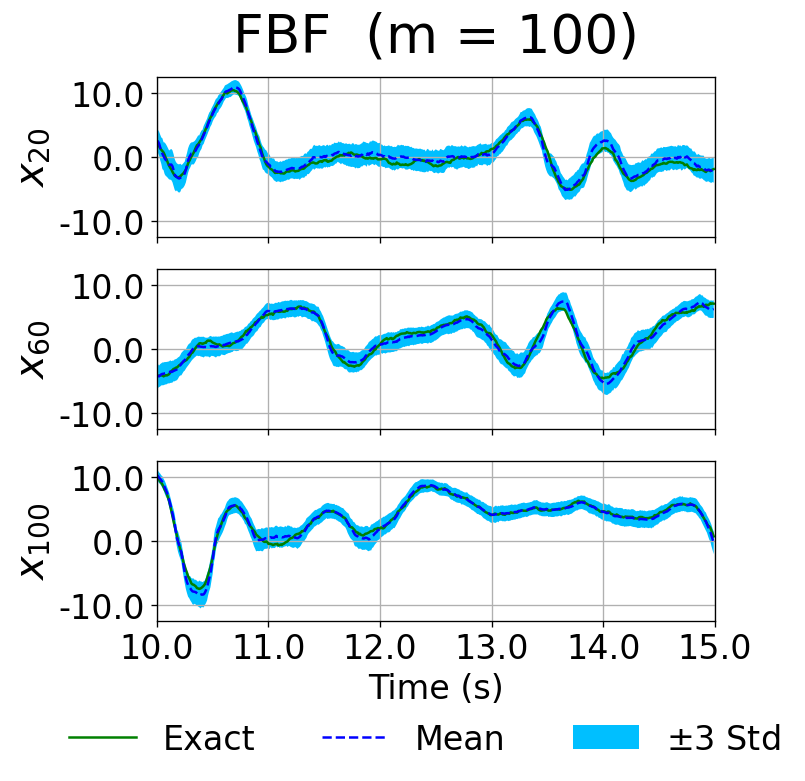}
	}
        \subfigure{
		\includegraphics[scale=0.42]{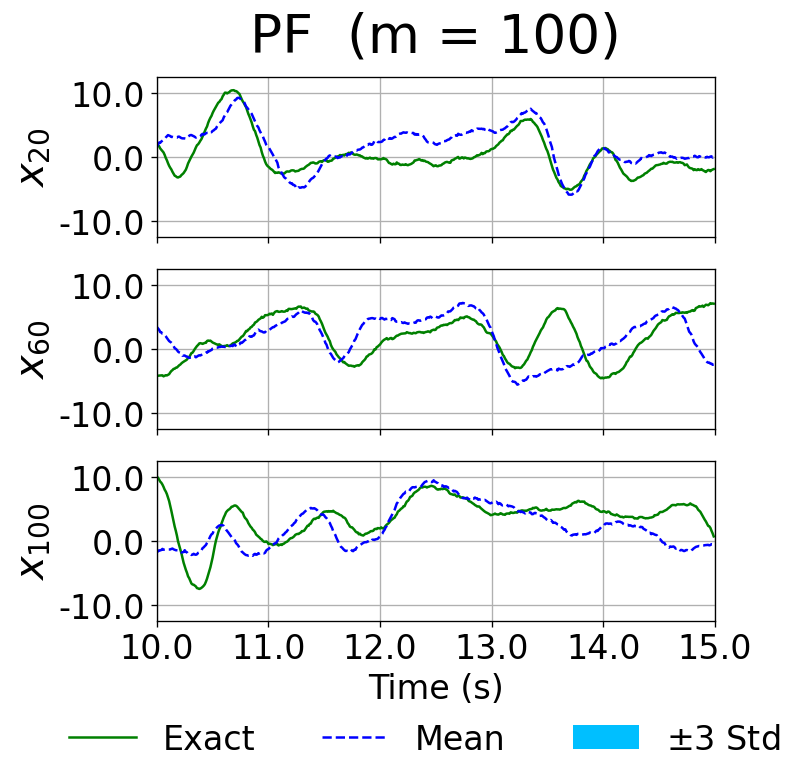}
	}
        \subfigure{
		\includegraphics[scale=0.42]{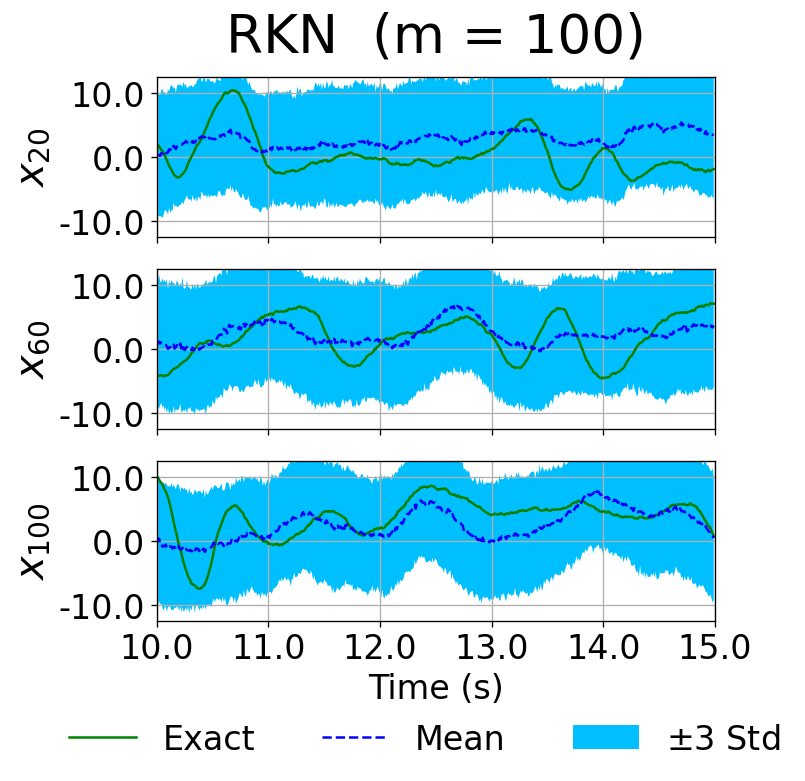}
	}
	\caption{Visualization of mean and uncertainty of the estimated posterior for a test measurement of $100$-dimensional Lorenz-96 systems. For clarity, we present the final 500 steps of the test trajectory of length $K_{\mathrm{test}} = 1500$.}
	\label{fig:ex2_post_interval_m=100}
\end{figure}

\begin{figure}[ht]
	\hspace{-20pt}
	\includegraphics[scale=0.46]{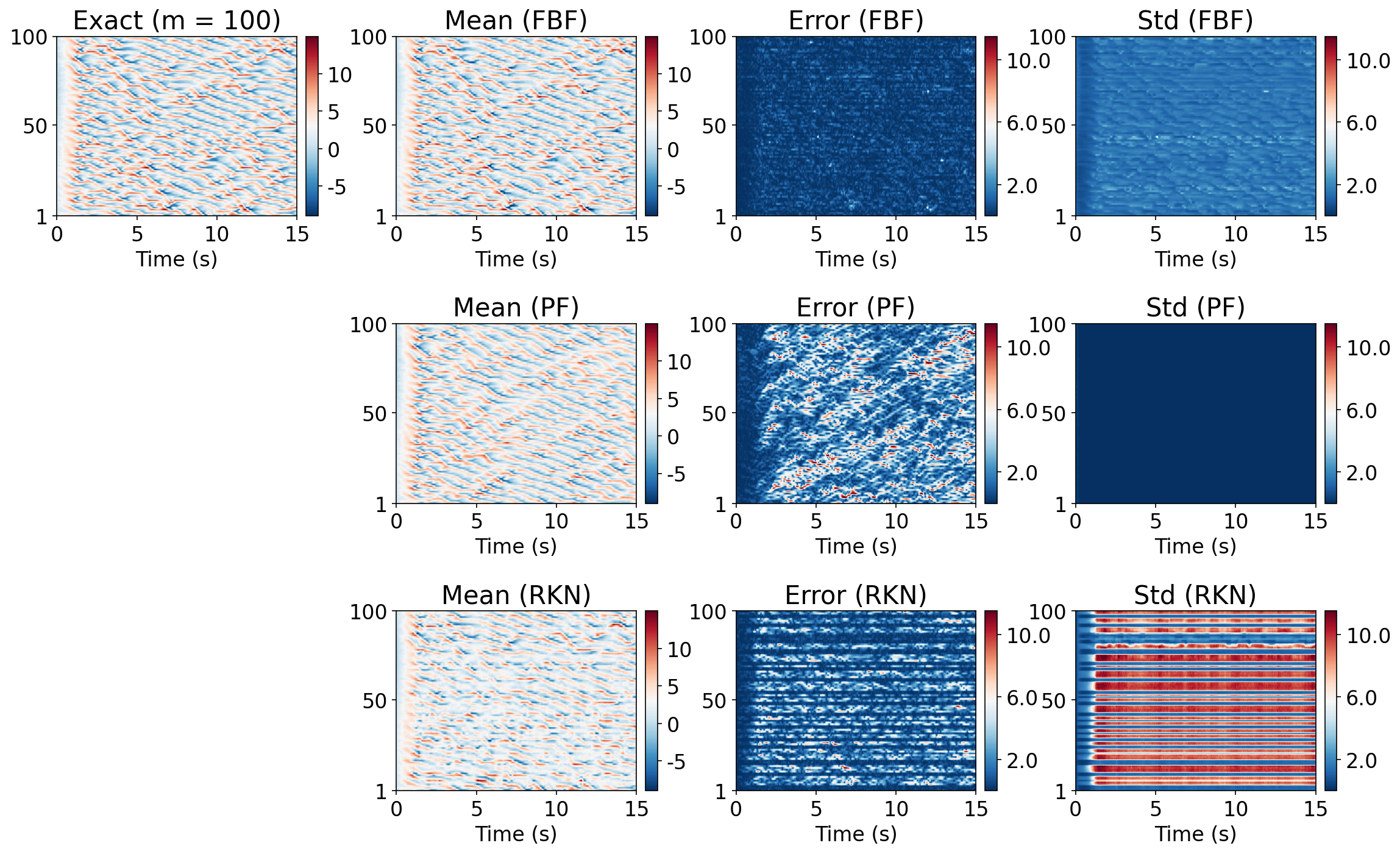}
	\caption{Results for state estimation of Lorenz-96 system with $m = 100$ on a test case: exact values (first column), mean and standard deviation of posterior samples (second and fourth columns), absolute errors between reference values and posterior mean (third column).}
	\label{fig:ex2_post_mesh_m=100}
\end{figure}

To visualize the effectiveness of high-dimensional state estimation using FBF, we present the mean and standard deviation of the posterior samples under two scenarios, with $m = 10$ and $m = 100$, in \cref{fig:ex2_post_interval_m=10} and \cref{fig:ex2_post_interval_m=100}. Note that \cref{fig:ex2_post_interval_m=10} and \cref{fig:ex2_post_interval_m=100} display three selected components while \cref{fig:ex2_post_mesh_m=100} portrays all components for the scenario $m = 100$. For the case of $m = 10$, we can see that (1) the CNF-DPF method produces significantly biased mean predictions; (2) the uncertainty estimates by the PF method fail to encompass the absolute error due to particle degeneration; (3) the mean predictions obtained by FBF and RKN closely approximate the exact solution, and the uncertainty estimates cover in most cases the point-wise errors. However, when $m$ increases to $100$, the RKN method exhibits significantly higher errors and uncertainties, and the PF method suffers from severe particle degeneracy. In contrast, the posterior mean obtained by FBF remains well aligned with the reference solution, and the standard deviation effectively captures the uncertainty.

Moreover, we present the computational cost for training both methods. It is evident that FBF requires significantly less offline time compared to the RKN method. This substantial difference can be attributed to the computational complexity of the recurrent Kalman filtering procedure employed by RKN during the offline stage. In contrast, the proposed FBF only necessitates learning a latent SSM, resulting in a more efficient training process.

\begin{table}[ht]
    \centering
    \caption{Comparison of the computational costs for different methods discussed in Section \ref{subsection:lorenz96_example}. The offline time corresponds to the training phase for each method, while the online time refers to the inference procedures for 200 test cases. For brevity, the offline and online times are denoted as "Off." and "On.", respectively. The least time (in seconds) among the methods is printed in bold font.}
    \begin{tabular}{ccccccccccccl}
    \hline
    \multicolumn{1}{l}{} & \multicolumn{2}{c}{FBF} &  & \multicolumn{2}{c}{PF} &  & \multicolumn{2}{c}{RKN} &  & \multicolumn{2}{c}{CNF-DPF} &  \\ \cline{2-3} \cline{5-6} \cline{8-9} \cline{11-12}
    \multicolumn{1}{l}{} & Off.            & On.   &  & Off.      & On.        &  & Off.     & On.          &  & Off.          & On.         &  \\ \hline
    $m = 10$             & $\bm{9419}$     & 590   &  & NA        & 6192       &  & 297987   & $\bm{357}$   &  & 55938         & 1861        &  \\
    $m = 20$             & $\bm{10927}$    & 602   &  & NA        & 11340      &  & 301422   & $\bm{361}$   &  & 59188         & 979         &  \\
    $m = 30$             & $\bm{12062}$    & 604   &  & NA        & 14270      &  & 309463   & $\bm{359}$   &  & 61895         & 755         &  \\
    $m = 40$             & $\bm{13115}$    & 791   &  & NA        & 17405      &  & 317216   & $\bm{385}$   &  & NA            & NA          &  \\
    $m = 50$             & $\bm{13113}$    & 850   &  & NA        & 20704      &  & 321176   & $\bm{394}$   &  & NA            & NA          &  \\
    $m = 100$            & $\bm{23704}$    & 968   &  & NA        &   $111081$         &  & 378972   &      $\bm{498}$        &  & NA            & NA          &  \\ \hline
    \end{tabular}
\end{table}

\subsection{Stochastic advection–diffusion model}\label{subsection:AD_example}
Finally, we experiment with the stochastic advection–diffusion equation driven by an additive random noise \cite{DEHGHAN:EABE:2015}, which has found extensive applications in hydrology \cite{Jinno:WRR:1993, Ancey:JGR:2015} and environmental monitoring \cite{DEHGHAN:EABE:2015}. Let us consider the following stochastic advection–diffusion equation:
\begin{equation}\label{equ:AD_example}
\left\{
\begin{array}{ll}
\mathrm{d}u + (\kappa\partial_{s}u - D\partial_{ss}^{2}u + g(s)){\mathrm{d}t} = \sigma\mathrm{d}W(t), \quad s\in [-1, 1],~~t\in [0, 1],\\
u(0, s) = -\mathrm{sin}(\pi s),\quad s\in [-1, 1],\\
u(t, -1) = 0,\quad t\in [0, 1],\\
u(t, 1) = 0,\quad t\in [0, 1],
\end{array}
\right.
\end{equation}
where $\kappa$ represents the advection coefficient, $D$ is the diffusion coefficient, $g(s)$ denotes the source function, $W(t)$ is a continuous Wiener process, and $u(s,t)$ indicates the concentration at position $s$ and time $t$. In this example, $g(s) = 5(s^2 - 1)$, $\sigma = 10$, and $D$ and $\kappa$ are set at $0.01$ and $0.5$, respectively. 

To generate the trajectory data, we simulate the advection-diffusion system using an implicit finite difference on a $201\times 100$ uniform temporal-spatial domain to obtain the reference solution for $u$. The target state and measurement are defined as:
\begin{eqnarray*}
\bm{x}_{k} & = & \left(u(k\Delta t,s_{1}),\ldots,u(k\Delta t,s_{100})\right),\\
\bm{y}_{k} & = & \left(e^{-u(k\Delta t,s_{1}^{\prime})-1},\ldots,e^{-u(k\Delta t,s_{n}^{\prime})-1}\right)+\bm{v}_{k},\quad\bm{v}_{k}\sim\mathcal{N}(0,r^{2}I_{n}),
\end{eqnarray*}
where $\Delta t = 0.005$, $\{s_1, \ldots, s_{100}\}$ are the spatial grid points in the implicit finite difference solution, and $\{s_1^{\prime}, \ldots, s_{n}^{\prime}\}$ are $n$ equally spaced sensor locations selected from $\{s_1, \ldots, s_{100}\}$ with $n < 100$. Here, we aim to quantify the uncertainty of the random field $u(k\delta t,\cdot)$ at $\{s_1, \ldots, s_{100}\}$ using filtering methods, based on incomplete and noisy observations $\bm{y}_k$ up to time $k\delta t$. We generate 1200 state and measurement pairs, dividing them into $1000$ trajectories for training and $200$ trajectories for testing. Additionally, the CNF-DPF is excluded from the comparison due to its excessive GPU memory requirement in high-dimensional situations.

\begin{figure}[!t]
	\centering
	\subfigure[]{
		\includegraphics[scale=0.4]{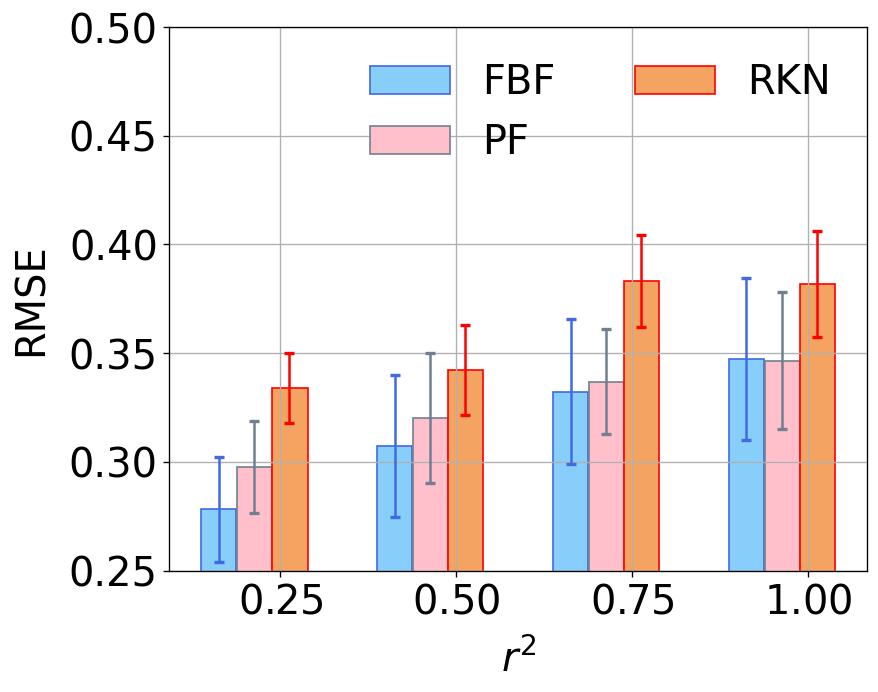}
	}
	\subfigure[]{
		\includegraphics[scale=0.4]{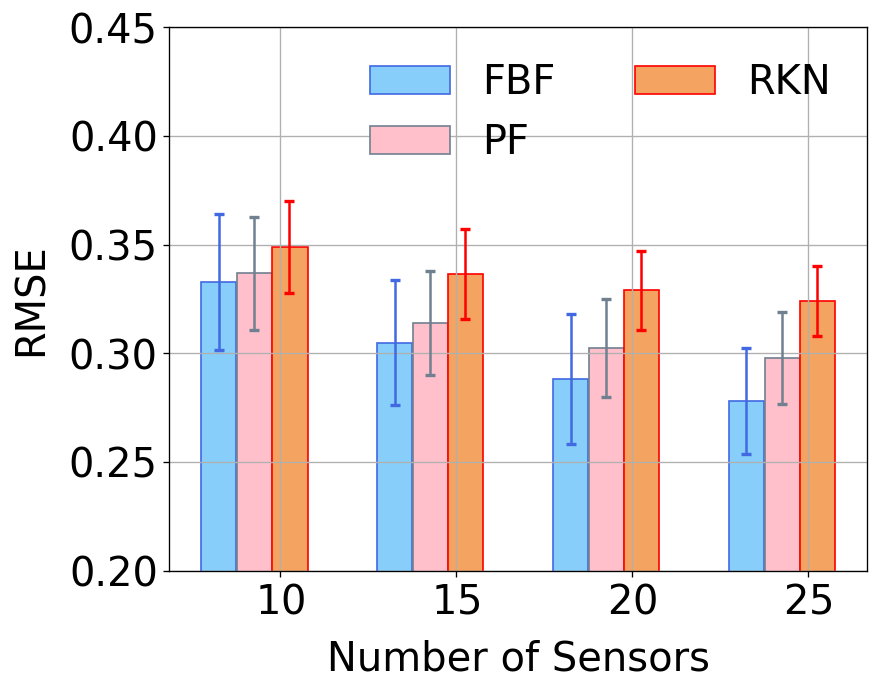}
	}
	\subfigure[]{
		\includegraphics[scale=0.4]{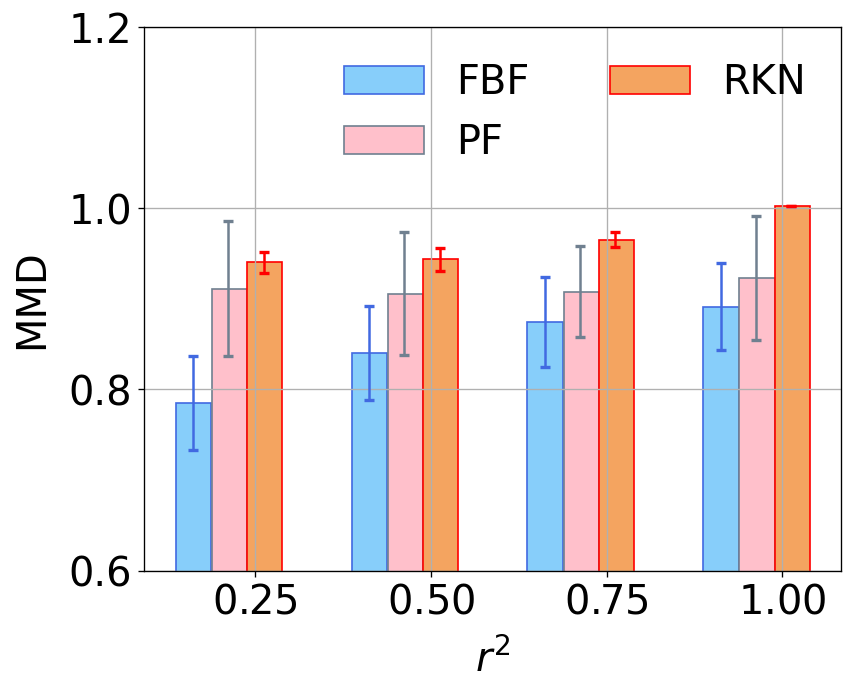}
	}
        \subfigure[]{
		\includegraphics[scale=0.4]{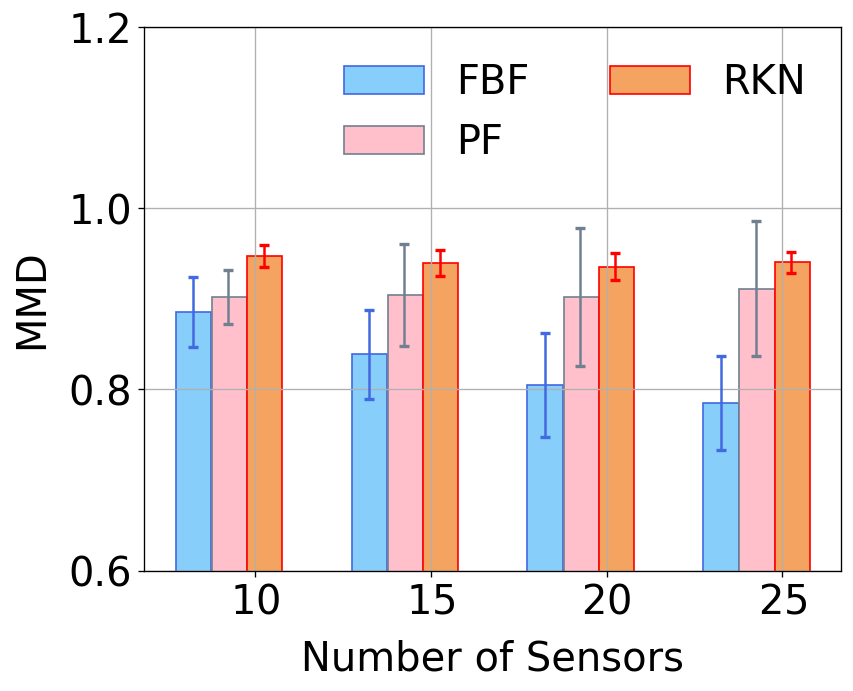}
	}
	\subfigure[]{
		\includegraphics[scale=0.4]{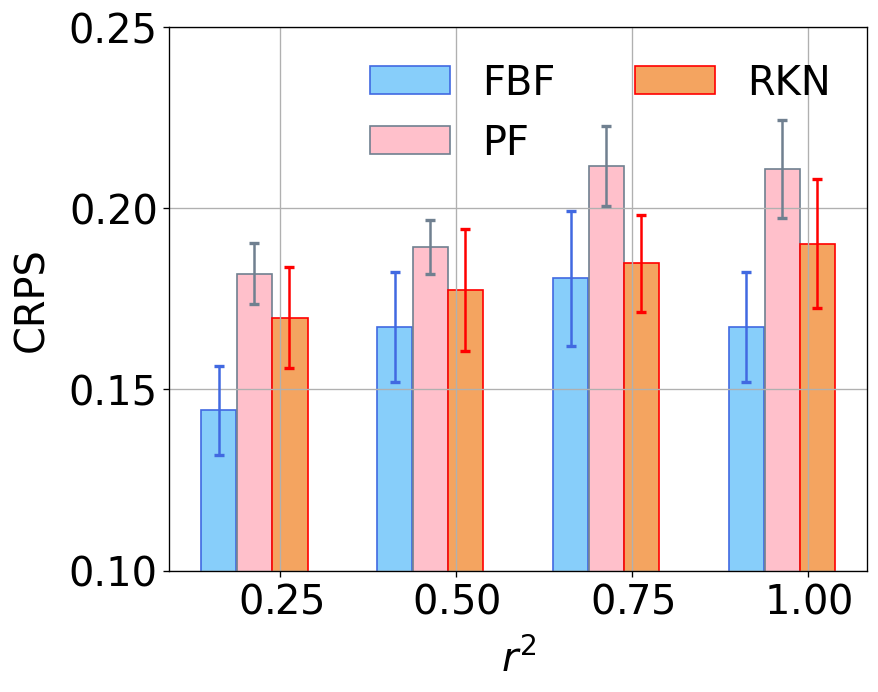}
	}
	\subfigure[]{
		\includegraphics[scale=0.4]{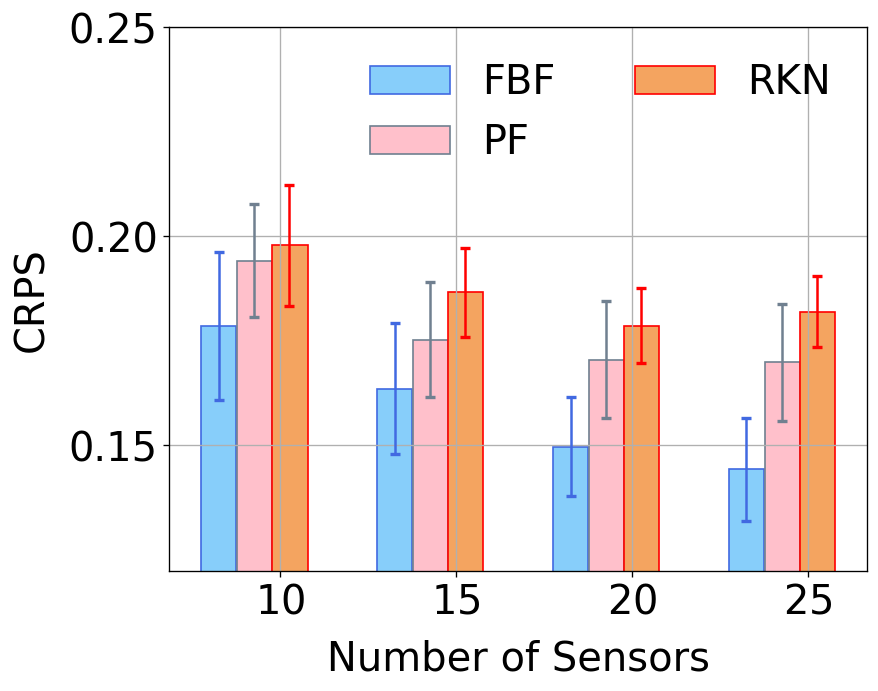}
	}
	\caption{Comparison of performance using different methods in \cref{subsection:AD_example}, where the bars represent the mean and the error bars indicate the standard deviation across 200 test cases. (a,b): RMSE; (c,d): MMD; (e,f): CRPS. Note that $r^2 = 0.25$ in (b,d,f).}
	\label{ex3_metric}
\end{figure}

We present the performance evaluation results for different methods with varying noise levels in \cref{ex3_metric} (a,c,e), where the number of sensors is fixed at $n = 25$. The results indicate that the FBF consistently achieves higher accuracy and provides more reliable UQ compared to other methods. To further validate the robustness of FBF, we reduce $n$ from $25$ to $10$ and present the evaluation results in \cref{ex3_metric} (b,d,f). While the performance of FBF declines with the reduction in observable components, it still outperforms the PF and RKN methods. 

\begin{figure}[!h]
	\centering
	\subfigure{
		\includegraphics[scale=0.58]{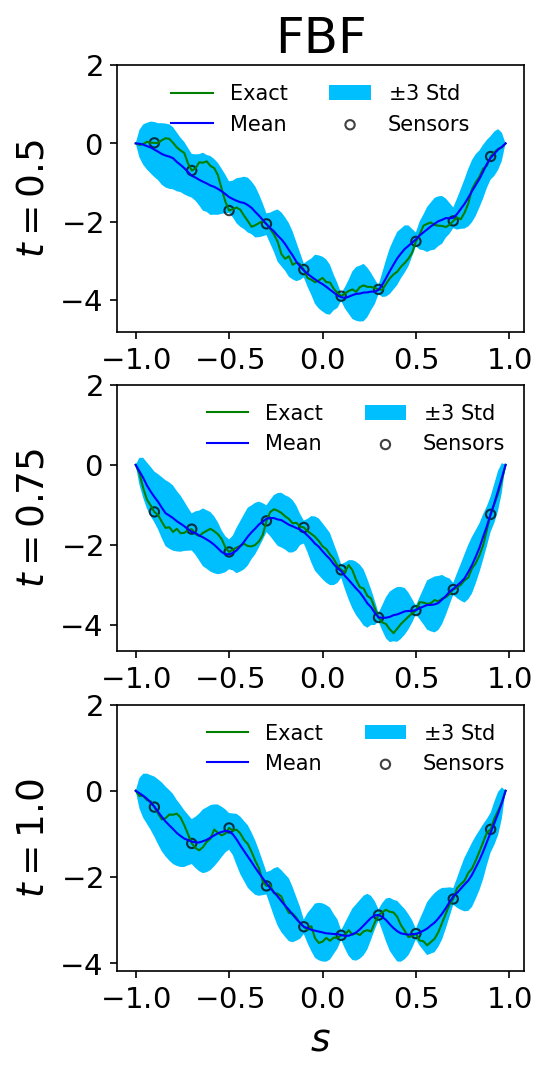}
	}
	\subfigure{
		\includegraphics[scale=0.58]{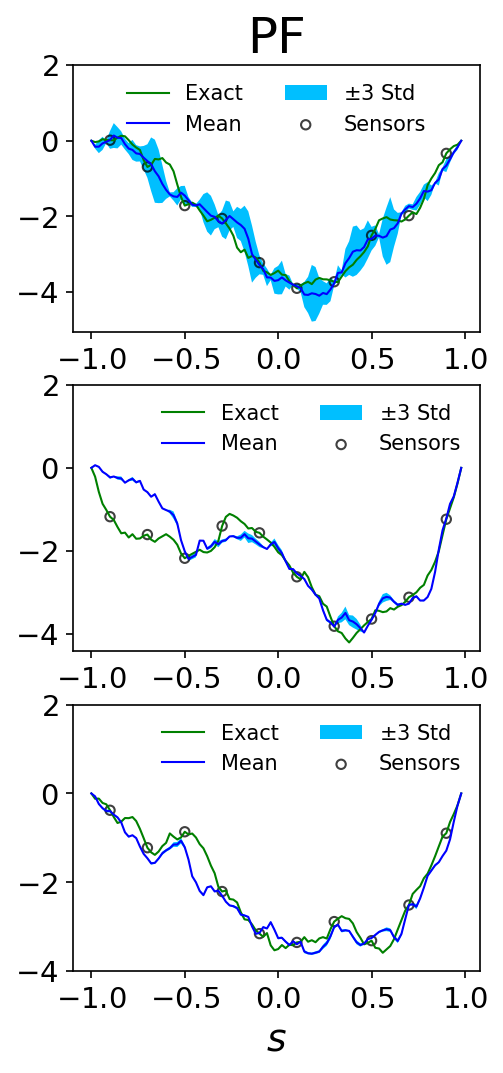}
	}
	\subfigure{
		\includegraphics[scale=0.58]{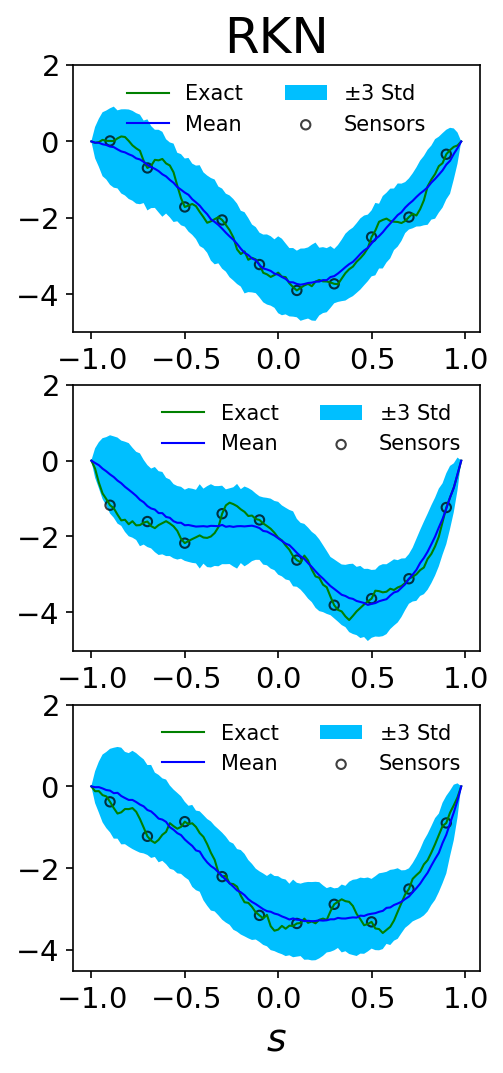}
	}
	\caption{Visualization of mean and uncertainty associated with the estimated posterior distributions obtained using three methods, given a test measurement consisting of 10 data points. The black circles indicate the sensor positions. Rows from top to bottom correspond to time $t = 0.5, 0.75, 1.0$, respectively. The number of samples utilized for PF is 10000, while 1000 samples are used in both FBF and RKN.}
	\label{ex3_CI}
\end{figure}

To visualize the effectiveness of UQ by the proposed FBF under sparse measurement conditions, \cref{ex3_CI} presents the mean predictions and uncertainty estimates obtained by the three methods for a test trajectory, with the number of sensors fixed at $n = 10$. It can be observed that FBF provides reliable uncertainty estimates for states. In contrast, the PF method suffers from a degeneracy issue, and its estimated uncertainty fails to cover the absolute errors, especially at time steps $t = 0.75$ and $t = 1.0$. While the RKN method provides uncertainty estimates that cover the absolute errors, it does not exhibit the expected reduction in standard deviation near sensor locations (indicated by black circles in \cref{ex3_CI}). In comparison, the FBF method offers more consistent performance, generating smaller standard deviations in proximity to sensor locations and larger standard deviations in regions without measurements. This highlights the capability of FBF to adapt uncertainty estimates based on measurement density, ensuring both accuracy and reliability.

Furthermore, the computational costs of the different methods are summarized in \cref{ex3_time_table}. It is observed that RKN necessitates expensive computational costs in the offline phase, while PF requires solving a $100$-dimensional linear system derived from an implicit difference scheme at each iteration for every particle, resulting in substantial computational burden during the online phase. In contrast, the proposed FBF achieves a balanced computational cost across the offline and online phases.

\begin{table}[!h]
	\centering
	\caption{Comparison of the computational cost required by different methods in Section \ref{subsection:AD_example}. The online time refers to the inference procedures for 200 test cases, each with 10 equally spaced sensors, while the offline time corresponds to 500 training iterations for FBF and RKN. The least time (in seconds) among the methods is printed in bold font.}
	\begin{tabular}{c|cccc}
		\hline
		Stage    & PF & FBF      & RKN     \\ \hline
		Offline  & NA & $\bm{895.26}$   & $95956.56$ \\
		Online   & $6443.25$ & $455.51$  & $\bm{162.12}$  \\ \hline
	\end{tabular}
	\label{ex3_time_table}
\end{table}

\FloatBarrier

\section{Conclusion}\label{section:conclusion}
In this paper, a novel flow-based Bayesian filter is proposed to efficiently solve nonlinear filtering problems. In contrast to existing deep learning techniques, the proposed method offers tractable estimation of the target filtering distribution, enabling efficient sampling and high-accuracy density evaluation. Numerical experiments demonstrate the advantages of our framework in terms of accuracy and efficiency, as well as its suitability for addressing high-dimensional filtering problems.
In future work, we will incorporate physical information into FBF to broaden its applicability and extend its use to real-world scenarios, including geophysical inversion and climate modeling.


\section*{Acknowledgement}
The authors were supported by the National Natural Science Foundation of China (No.~12271409), the Foundation of National Key Laboratory of Computational Physics, and the Fundamental Research Funds for the Central Universities.

\bibliography{refs}

\appendix
\section{A variant of FBF}\label{appendix:FBF-cSSM} 
\subsection{Model structure}
This section introduces a variant of the FBF, referred to as FBF', where the latent SSM is formulated as a linear SSM with Gaussian noise, a commonly used model in filtering studies. The latent SSM in FBF' is expressed as:
\begin{eqnarray}
\bm{\chi}_{k} &=& E + F\bm{\chi}_{k-1} + \bm{\zeta}_{k}, \quad \bm{\zeta}_{k} \sim \mathcal{N}(0, P_{\chi}), \label{eq:latent-update-cSSM} \\
\bm{\gamma}_k &=& G + H\bm{\chi}_{k} + \bm{\xi}_{k}, \quad \bm{\xi}_{k} \sim \mathcal{N}(0, P_{\gamma}), \label{eq:latent-prediction-cSSM}
\end{eqnarray}
Here, \( F \in \mathbb{R}^{m \times m} \) and \( H \in \mathbb{R}^{n \times m} \) are the transition and measurement matrices, \( E \in \mathbb{R}^{m} \) and \( G \in \mathbb{R}^{n} \) are constant drift terms, and \( \bm{\zeta}_{k} \) and \( \bm{\xi}_{k} \) represent independent Gaussian noise with diagonal covariance matrices \( P_{\chi} \) and \( P_{\gamma} \), respectively.
\( E, F, G \) and \( H \) are all trainable and independent of the state and observation variables.


\subsection{Training and filtering algorithm}
According to the latent SSM in \eqref{eq:latent-update-cSSM} and \eqref{eq:latent-prediction-cSSM}, the conditional probabilities $p(\bm{x}_{k}|\bm{x}_{k-1})$ and $p(\bm{y}_{k}|\bm{x}_{k})$ can be evaluated as follows:
\begin{eqnarray}
p(\bm{x}_{k}|\bm{x}_{k-1}) & = & p(\bm{\chi}_{k}|\bm{\chi}_{k-1})\left|\mathrm{det}\left(\frac{\partial\mathcal{T}(\bm{x}_{k})}{\partial\bm{x}_{k}}\right)\right|\nonumber \\
 & = & \mathcal{N}\left(\bm{\chi}_{k}|E + F\bm{\chi}_{k-1},P_{\gamma}\right)\left|\mathrm{det}\left(\frac{\partial\mathcal{T}(\bm{x}_{k})}{\partial\bm{x}_{k}}\right)\right|\nonumber \\
 & = & \mathcal{N}\left(\mathcal{T}(\bm{x}_k)|E + F\mathcal{T}(\bm{x}_{k-1}),P_{\gamma}\right)\left|\mathrm{det}\left(\frac{\partial\mathcal{T}(\bm{x}_{k})}{\partial\bm{x}_{k}}\right)\right|\nonumber \\
 & \triangleq & f_{s}^{'}(\bm{x}_{k-1},\bm{x}_{k}),\\
p(\bm{y}_{k}|\bm{x}_{k}) & = & p(\bm{\gamma}_{k}|\bm{\chi}_{k})\left|\mathrm{det}\left(\frac{\partial\mathcal{V}(\bm{y}_{k})}{\partial\bm{y}_{k}}\right)\right|\nonumber \\
 & = & \mathcal{N}\left(\bm{\gamma}_{k}|G + H\bm{\chi}_{k},P_{\gamma}\right)\left|\mathrm{det}\left(\frac{\partial\mathcal{V}(\bm{y}_{k})}{\partial\bm{y}_{k}}\right)\right|\nonumber \\
  & = & \mathcal{N}\left(\mathcal{V}(\bm{y}_{k})|G + H\mathcal{T}(\bm{x}_{k}),P_{\gamma}\right)\left|\mathrm{det}\left(\frac{\partial\mathcal{V}(\bm{y}_{k})}{\partial\bm{y}_{k}}\right)\right|\nonumber \\
 & \triangleq & f_{o}^{'}(\bm{x}_{k},\bm{y}_{k}).
\end{eqnarray}
Suppose we are given a data set \( \{\bm{x}_{0:K}^{\mathrm{train}}, \bm{y}_{1:K}^{\mathrm{train}}\} \). Similar to FBF, all model parameters in FBF' can be optimized by maximizing the following objective function:
\[
\mathcal{L} = \frac{\alpha^{'}}{K} \sum_{k=1}^{K} \log f_s^{'}(\bm{x}_{k-1}^{\mathrm{train}}, \bm{x}_k^{\mathrm{train}}) + \frac{\beta^{'}}{K} \sum_{k=1}^{K} \log f_o^{'}(\bm{x}_{k}^{\mathrm{train}}, \bm{y}_k^{\mathrm{train}}),
\]
where \( \alpha^{'} \) and \( \beta^{'} \) are weight factors. The training procedure is analogous to that of FBF and is therefore omitted here for brevity.

In the filtering procedure, if the latent filtering distribution \( p(\bm{\chi}_{k-1}~|~\bm{\gamma}_{1:k-1}) \) at time step \( k-1 \) is given as a Gaussian distribution with mean \( \bm{\mu}_{k-1} \) and covariance \( \bm{\Sigma}_{k-1} \), the latent filtering distribution at time step \( k \) can be computed using the Kalman filter \cite{Kalman:JBE:1960} as follows, when a new measurement \( \bm{\gamma}_k \) becomes available:
\begin{equation}\label{equ:KF}
\begin{aligned}
\bm{\mu}_{k} &= E + F\bm{\mu}_{k-1} + K_k(\bm{\gamma}_k - HF\bm{\mu}_{k-1} - HE - G),\\
\bm{\Sigma}_{k} &= (I_m - K_kH)(F\bm{\Sigma}_{k-1}F^{\top} + P_{\chi}),
\end{aligned}
\end{equation}
where the Kalman gain \( K_k \) and the innovation covariance \( S_k \) are defined as:
\[
\begin{aligned}
K_k &= (F\bm{\Sigma}_{k-1}F^{\top} + P_{\chi})H^{\top}S_k^{-1},\\
S_k &= H(F\bm{\Sigma}_{k-1}F^{\top} + P_{\chi})H^{\top} + P_{\gamma}.
\end{aligned}
\]
Subsequently, samples from the target filtering distribution $p(\bm{x}_k~|~\bm{y}_{1:k})$ can be produced by \eqref{equ:FBF_infer}.

\subsection{Comparison with FBF}\label{appendix:proof-SSM}
According to \eqref{eq:latent-update-cSSM} and \eqref{eq:latent-prediction-cSSM}, the conditional probabilities $p(\bm{\chi}_{k}~|~\bm{\chi}_{k-1},\bm{\gamma}_{k-1})$ and $p(\bm{\gamma}_{k}~|~\bm{\chi}_{k-1})$ provided by the FBF' can be formulated as:
\begin{align*}
p(\bm{\chi}_{k}~|~\bm{\chi}_{k-1},\bm{\gamma}_{k}) & \propto p(\bm{\gamma}_{k}~|~\bm{\chi}_{k})\cdot p(\bm{\chi}_{k}~|~\bm{\chi}_{k-1})\\
 & =\mathcal{N}(\bm{\gamma}_{k}~|~G+H\bm{\chi}_{k},P_{\gamma})\cdot\mathcal{N}(\bm{\chi}_{k}~|~E+F\bm{\chi}_{k-1},P_{\chi})\\
 & \propto\mathcal{N}(\bm{\gamma}_{k}~|~\bm{\mu}_{k}^{\prime\prime},\bm{\Sigma}_{k}^{\prime\prime}),\\
p(\bm{\gamma}_{k}~|~\bm{\chi}_{k-1}) & =\int p(\bm{\gamma}_{k}~|~\bm{\chi}_{k})~p(\bm{\chi}_{k}~|~\bm{\chi}_{k-1})\mathrm{d}\bm{\chi}_{k}\\
 & =\int\mathcal{N}(\bm{\gamma}_{k}~|~G+H\bm{\chi}_{k},P_{\gamma})~\mathcal{N}(\bm{\chi}_{k}~|~E+F\bm{\chi}_{k-1},P_{\chi})\mathrm{d}\bm{\chi}_{k}\\
 & =\mathcal{N}(\bm{\gamma}_{k}~|~G+HE+HF\bm{\chi}_{k-1},P_{\gamma}+HP_{\chi}H^{\top}),
\end{align*}
where
\begin{align*}
\bm{\mu}_{k}^{\prime\prime} & =E-K_{k}^{'}(G+HE)+K_{k}^{'}\bm{\gamma}_{k}+(F-K_{k}^{'}HF)\bm{\chi}_{k-1},\\
\bm{\Sigma}_{k}^{\prime\prime} & =(I_{m}-K_{k}^{'}H)P_{\chi},\\
\ensuremath{K_{k}^{'}} & =P_{\chi}H^{\top}(HP_{\chi}H^{\top}+P_{\gamma})^{-1}.
\end{align*}
Thus, the latent SSM in FBF' can be equivalently expressed in the form of \eqref{eq:latent-update} and \eqref{eq:latent-prediction}:
\begin{align}
\bm{\chi}_{k} &= \underbrace{E - K_{k}^{'}(G+HE) + K_k^{'}\bm\gamma_{k}}_{A} + \underbrace{(F - K_{k}^{'}HF)}_{B}\bm\chi_{k-1} + \bm\zeta_{k}^{'}, \quad \bm\zeta_{k}^{'}\sim\mathcal{N}(0, \underbrace{(I_m - K_k^{'}H)P_{\chi}}_{Q_{\chi}}),\label{latent-update-cSSM-equiv}\\
\bm{\gamma}_{k} &= \underbrace{G+HE}_{C} + \underbrace{HF}_{D}\bm{\chi}_{k-1} + \bm{\xi}_{k}^{'}, \quad\bm{\xi}_{k}^{'}\sim\mathcal{N}(0, \underbrace{HP_{\chi}H^{\top} + P_{\gamma}}_{Q_{\gamma}}).\label{latent-prediction-cSSM-equiv}
\end{align}

From the above analysis, it can be seen that while both latent SSMs in FBF and FBF' preserve the Gaussianity of the latent filtering distribution, the latent SSM in FBF is more general than that in FBF'. Additionally, our numerical experiments also indicate that FBF achieves better filtering performance (see below).

\subsection{Numerical Results}
We apply FBF' to the filtering problem in \cref{subsection:sinusoidal-example} and \cref{subsection:lorenz96_example}, and compare its performance with FBF. The performance evaluation is provided in \cref{Appendix_test_metric}, and the computational cost is detailed in \cref{Appendix_time_table}.

We can observe that the computational cost of FBF$^{'}$ remains on the same order of magnitude as FBF, but FBF$^{'}$ exhibits higher errors than FBF. This demonstrates the advantages brought by the generality of our proposed latent SSM compared to the traditional linear SSM.

\begin{figure}[!h]
	\centering
	\subfigure[]{
		\includegraphics[scale=0.4]{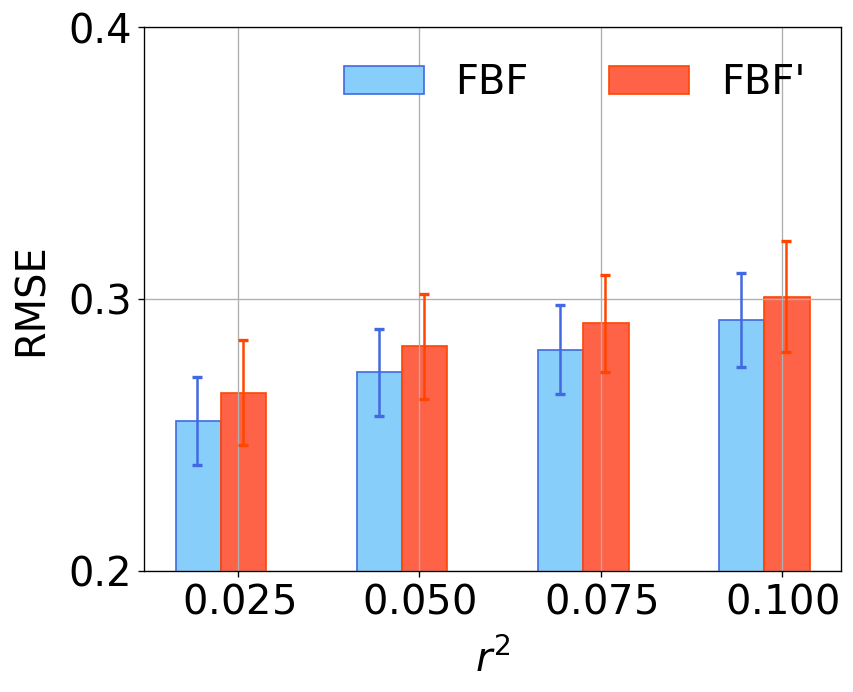}
	}
	\subfigure[]{
		\includegraphics[scale=0.4]{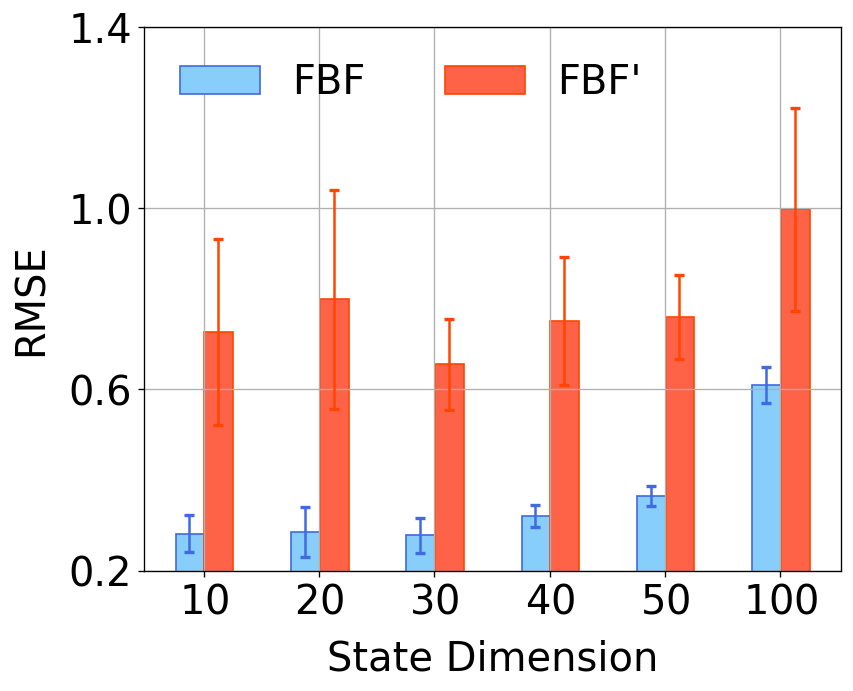}
	}
	\subfigure[]{
		\includegraphics[scale=0.4]{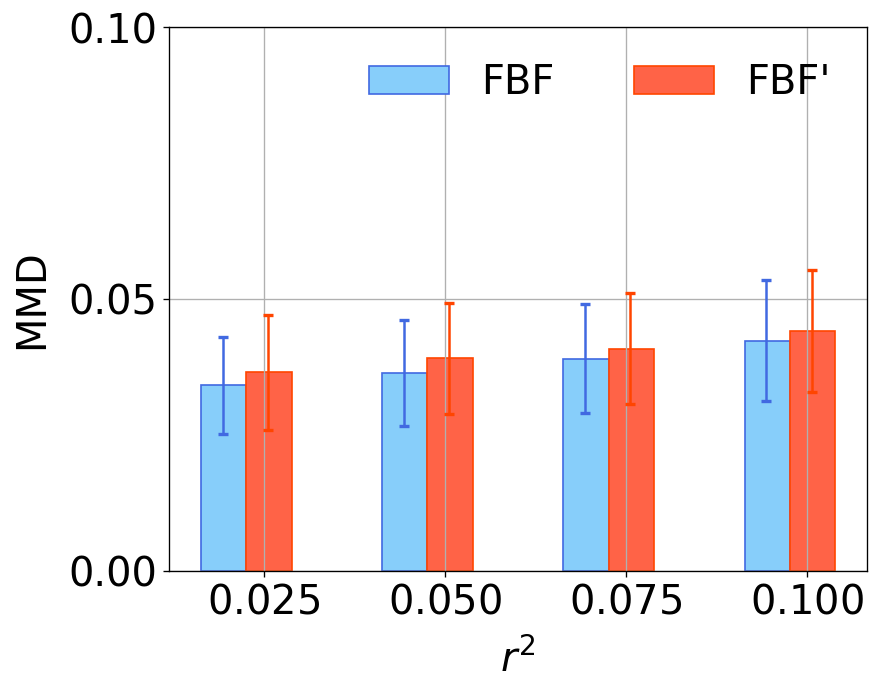}
	}
        \subfigure[]{
		\includegraphics[scale=0.4]{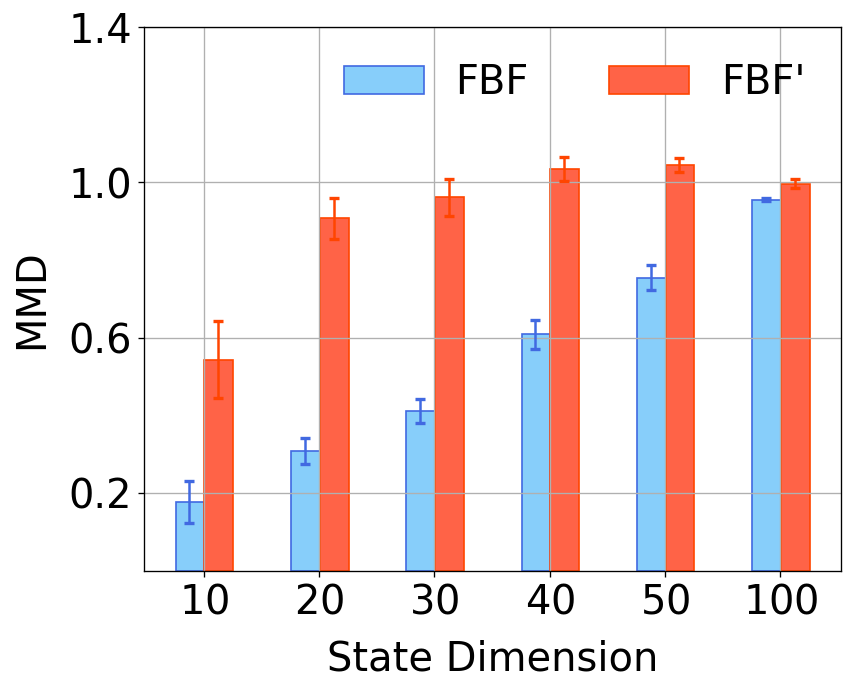}
	}
	\subfigure[]{
		\includegraphics[scale=0.4]{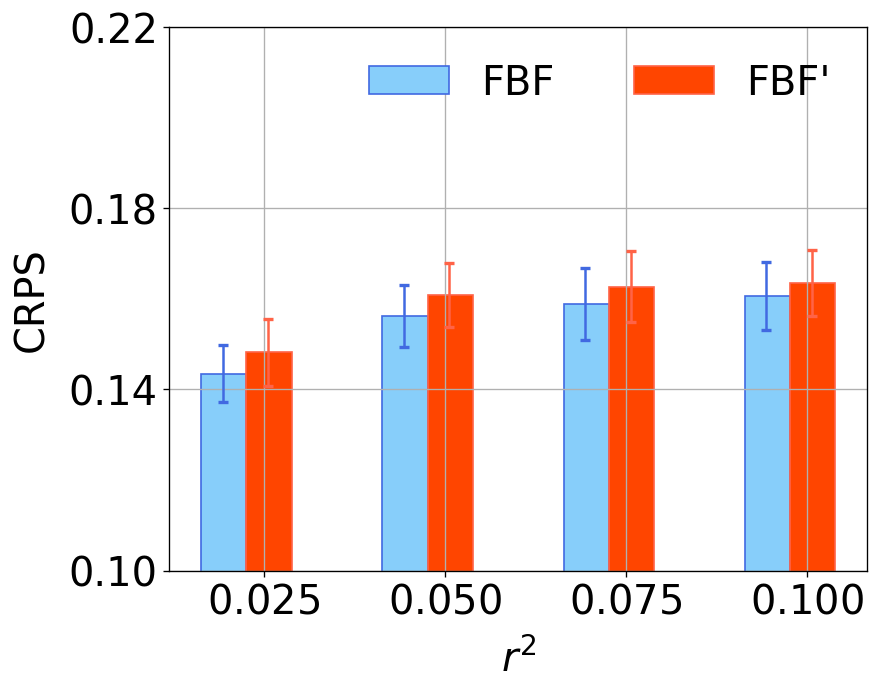}
	}
	\subfigure[]{
		\includegraphics[scale=0.4]{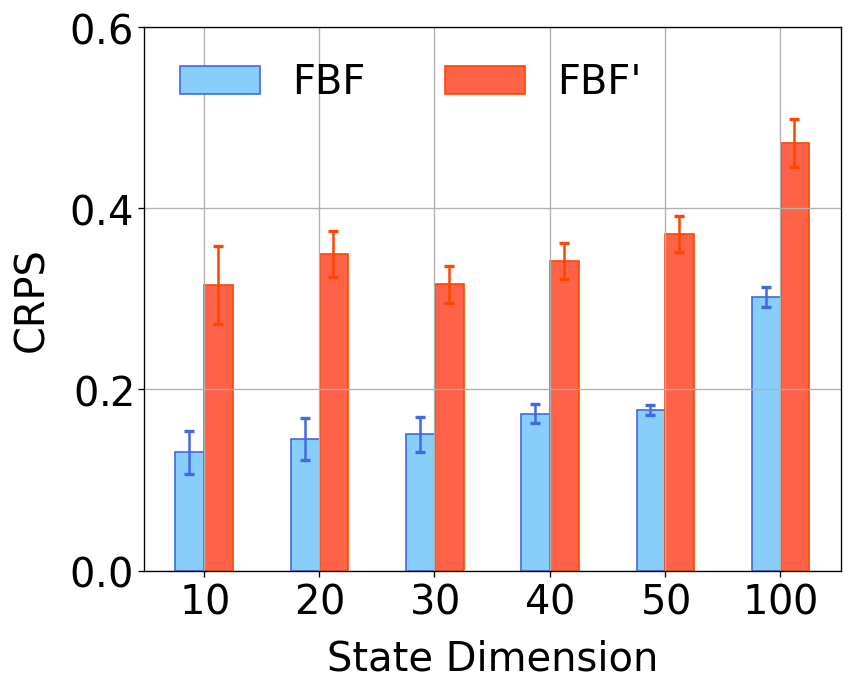}
	}
	\caption{Comparison of performance using FBF and FBF$^{'}$. (a),(c),(e): metric values of experiments in \cref{subsection:sinusoidal-example}; (b),(d),(f): metric values of experiments in \cref{subsection:lorenz96_example}.}
	\label{Appendix_test_metric}
\end{figure}

\begin{table}[!h]
	\centering
	\caption{Comparison of the computational cost required by FBF$^{'}$ and FBF. The offline time reflects $500$ training iterations, and the online time represents the inference procedures for test cases.}
	\begin{tabular}{c|ccccc}
		\hline
		Stage    & FBF      & FBF$^{'}$     \\ \hline
		Offline  & $435.18$ & $469.97$ \\
		Online   & $120.44$ & $132.81$ \\ \hline
	\end{tabular}
	\label{Appendix_time_table}
\end{table}

\section{Details for computing all the metrics}\label{appendix:metrics}
\subsection{RMSE}
The root mean square error (RMSE) between the mean values of the estimated filtering distributions and the true state values, denoted as $\bm{x}_k^{\mathrm{true}}$, is defined as \cite{HJ:2018}:
\[
\mathrm{RMSE} = \sqrt{\frac{1}{mK} \sum_{k=1}^{K} \left\Vert \bm{x}_k^{\mathrm{true}} - \frac{1}{N} \sum_{j=1}^{N} \bm{x}_k^{(j)} \right\Vert^2},
\]
where $\Vert\cdot\Vert$ represents the Euclidean norm of vectors, $m$ is the state dimensionality, and $K$ represents the total number of time steps.

\subsection{MMD}
Let $\phi(\cdot)$ denote a feature mapping derived from a kernel function $\mathrm{ker}(\cdot,\cdot)$ such that $\phi(\bm{x})^\top \phi(\bm{y}) = \mathrm{ker}(\bm{x}, \bm{y})$. At each time step $k$, the maximum mean discrepancy (MMD) \cite{Gretton:JMLR:2012} between the empirical distribution of samples $\{\bm{x}_k^{(j)}\}_{j=1}^N$ generated by the filter and the delta distribution centered at $\bm{x}_k^{\mathrm{true}}$ is given by
\begin{eqnarray*}
    \mathrm{MMD}_k &=& \left\Vert \phi(\bm{x}_k^{\mathrm{true}}) - \frac{1}{N}\sum_{j=1}^N \phi(\bm{x}_k^{(j)}) \right\Vert^2 \\
    &=& \frac{1}{N^2} \sum_{i,j} \mathrm{ker}(\bm{x}_k^{(i)}, \bm{x}_k^{(j)}) - \frac{2}{N} \sum_{j} \mathrm{ker}(\bm{x}_k^{(j)}, \bm{x}_k^{\mathrm{true}})+\mathrm{ker}(\bm x_k^{\mathrm{true}},\bm x_k^{\mathrm{true}}),
\end{eqnarray*}
which quantifies the discrepancy between the expectation of the estimated filtering distribution and the true state in the feature space defined by $\phi$.

To compare the performance of different filters, we use the average MMD:
\[
\mathrm{MMD} = \frac{1}{K} \sum_{k=1}^K \mathrm{MMD}_k.
\]

In this paper, we select the kernel function as $\mathrm{ker}(\bm{x},\bm{y}) = \mathrm{exp}\Big(-\dfrac{\Vert\bm{x} - \bm{y}\Vert^2}{2\sigma^2}\Big)$, where $\sigma = 2$ for all experiments.


\subsection{CRPS}
The continuous ranked probability score (CRPS) \cite{Jochen:JRMS:2012} is commonly used to measure the difference between the cumulative distribution of predicted samples and the empirical distribution of the true state in the context of state estimation. For the $i$-th component of the state at time step $k$, it is defined as
\[
\mathrm{CRPS}_{k,i} = \int_{-\infty}^{+\infty}\left( \mathbf{1}(x_{k}^{\mathrm{true}} < x) - \frac{1}{N}\sum_{j=1}^N \mathbf{1}(x_{k,i}^{(j)} < x) \right)^2 \, \mathrm{d}x.
\]
Here $x_{k,i}^{(j)}$ denotes the $i$-th component of the $j$-th sample $\bm{x}_k^{(j)}$, and $\mathbf{1}(\cdot)$ is the indicator function, which equals 1 if the condition inside the parentheses is satisfied, and 0 otherwise.  

To assess the overall performance of filter methods, we compute the CRPS by averaging $\mathrm{CRPS}_{k,i}$ across state components and time steps \cite{SBM:SR:2022}:
\[
\mathrm{CRPS} = \dfrac{1}{mK}\sum_{k=1}^{K}\sum_{i=1}^{m}\mathrm{CRPS}_{k,i}.
\]

\section{Details for numerical implementations}\label{appendix:hyperparams}
\subsection{Calculation of covariance matrices}
 For ease of practical implementation, we assume $Q_{\chi}(\bm\gamma_k),Q_{\gamma}$ in \eqref{eq:latent-update} and \eqref{eq:latent-prediction}, as well as $P_{\chi},P_{\gamma}$ in \eqref{eq:latent-update-cSSM} and \eqref{eq:latent-prediction-cSSM}, are parameterized as diagonal covariance matrices. Specifically, $Q_{\chi}(\bm\gamma_k)$ is calculated as
 \[
Q_{\chi}(\bm\gamma_k) = \begin{bmatrix}
 \mathrm{Softplus}(\tilde{Q}_{\chi,1}) & &\\
 & \ddots & \\
 & & \mathrm{Softplus}(\tilde{Q}_{\chi,m})
\end{bmatrix},
 \]
where $\mathrm{Softplus}$ function ensures the positive definiteness of covariance matrices, and the terms $(\tilde{Q}_{\chi,1},\cdots,\tilde{Q}_{\chi,m})$ denotes the output vector of neural network $\tilde{Q}(\bm{\gamma}_k)$, whereas for $Q_{\gamma}, P_{\chi}, P_{\gamma}$, their diagonal elements are set to trainable parameters. Although a full covariance representation offers greater generality, our numerical experiments demonstrate that this diagonal parameterization does not affect the performance of our method.

\subsection{Neural networks architecture and training hyperparameters}
Subsequently, we outline the hyperparameter settings for all the methods employed in the numerical examples. For the PF methods used in \cref{subsection:sinusoidal-example}, \cref{subsection:lorenz96_example} and \cref{subsection:AD_example}, the number of particles is set to $2000$, $20000$, and $10000$, respectively. In contrast, the hyperparameters of network architecture for all data-driven methods are described among \cref{table:FBF} to \cref{table:DPF_ex1_ex2}, where \cref{table:FBF} presents the hyperparameters of FBF, \cref{table:RKN_ex1_ex2} and \cref{table:RKN_ex3} list the hyperparameters of RKN, and the hyperparameters of CNF-DPF are provided in \cref{table:DPF_ex1_ex2}.

Additionally, we present the training hyparameters in \cref{table:training}. Note that the weighting factors for the proposed FBF method are set to $1$, i.e., $\alpha = \beta = 1$, in all numerical experiments. During the training process, the learning rate for all methods decays according to the following formula \cite{ardizzone2018analyzing}:
\[
\eta_{\tau} = \eta_0\cdot(d)^{\frac{\tau}{E_{\text{train}}}},
\]
where $\eta_{\tau}$ denotes the learning rate at $\tau\mathrm{-th}$ training iterations, $\eta_0$ represents the initial learning rate, and $d$ is the final decaying factor, which is set to $0.1$ for all the methods.

\vspace{+5pt}
\noindent\textbf{Abbreviations}: $N_{\text{block}}$, the number of affine coupling blocks; $N_{\text{layer}}$, the number of fully-connected layers; $N_{\text{units}}$, the number of hidden units each layer of MLP; Act., activation functions; $E_{\text{epoch}}$: training epochs; $N_{\text{batch}}$: batch size; lr: initial learning rate.

\begin{table}[h]
	\centering
	\caption{Hyperparameters of the network structure for the FBF method.}
	\begin{tabular}{lccccccc}
		\hline
		& \multicolumn{3}{c}{NF} & & \multicolumn{3}{c}{DNN} \\ \cline{2-4} \cline{6-8} 
		& $N_{\text{block}}$ & $N_{\text{layer}}$ & $N_{\text{units}}$ & & $N_{\text{layer}}$ & $N_{\text{units}}$ & Act. \\ \hline
		\cref{subsection:sinusoidal-example} & 6 & 3 & 64 & & 6 & 64 & ReLU \\
		\cref{subsection:lorenz96_example}   & 8 & 3 & 128 & & 8 & 128 & ReLU \\
		\cref{subsection:AD_example}         & 6 & 3 & 64 & & 6 & 64 & ReLU \\ \hline
	\end{tabular}
	\label{table:FBF}
\end{table}

\begin{table}[!h]
	\centering
	\caption{Hyperparameters for the RKN method in \cref{subsection:sinusoidal-example} and \cref{subsection:AD_example}. Notations: $N_{\text{basis}}$: the number of basis for constructing transition matrices; bandwidth: bandwidth for band matrices; $N_{\text{lod}}$: the dimensionality of latent observations.}
	\begin{tabular}{cccclccclccc}
		\hline
		\multicolumn{1}{l}{} & \multicolumn{3}{c}{Encoder}                     &  & \multicolumn{3}{c}{Decoder}                     &  & \multicolumn{3}{c}{RKN Layer}                     \\ \cline{2-4} \cline{6-8} \cline{10-12} 
		& $N_{\text{layer}}$ & $N_{\text{units}}$ & Act. &  & $N_{\text{layer}}$ & $N_{\text{units}}$ & Act. &  & $N_{\text{basis}}$ & bandwidth & $N_{\text{lod}}$ \\ \hline
		\cref{subsection:sinusoidal-example}              & 4                  & 64                  & ReLU &  & 3                  & 64                  & Tanh &  & 8                  & 2         & 2                \\
		\cref{subsection:AD_example}              & 4                  & 64                  & ReLU & & 3                  & 64                  & Tanh &  & 8                  & 3         & 3                \\ \hline
	\end{tabular}
	\label{table:RKN_ex1_ex2}
\end{table}

\begin{table}[!h]
	\centering
	\caption{Hyperparameter of the network structure for the RKN method under varying state dimensionalities in \cref{subsection:lorenz96_example}.}
	\begin{tabular}{lccccccclccc}
		\hline
		& \multicolumn{3}{c}{Encoder}                     & \multicolumn{1}{l}{} & \multicolumn{3}{c}{Decoder}                                         &  & \multicolumn{3}{c}{RKN Layer}                     \\ \cline{2-4} \cline{6-8} \cline{10-12} 
		\multicolumn{1}{c}{}         & $N_{\text{layer}}$ & $N_{\text{units}}$ & Act. & \multicolumn{1}{l}{} & $N_{\text{layer}}$ & $N_{\text{units}}$ & Act.                     &  & $N_{\text{basis}}$ & bandwidth & $N_{\text{lod}}$ \\ \hline
		\multicolumn{1}{c}{$m = 10$} & 5                  & 64                  & ReLU &                      & 3                  & 64                  & Tanh                     &  & 8                  & 6         & 10               \\
		\multicolumn{1}{c}{$m = 20$} & 5                  & 64                  & ReLU &                      & 3                  & 64                  & Tanh                     &  & 8                  & 15        & 20               \\
		$m = 30$                     & 5                  & 100                 & ReLU &                      & 3                  & 100                 & \multicolumn{1}{l}{Tanh} &  & 8                  & 15        & 30               \\
		$m = 40$                     & 5                  & 100                 & ReLU &                      & 3                  & 100                 & \multicolumn{1}{l}{Tanh} &  & 8                  & 25        & 40               \\
		$m = 50$                     & 5                  & 100                 & ReLU &                      & 3                  & 100                 & \multicolumn{1}{l}{Tanh} &  & 8                  & 30        & 50               \\
        $m = 100$                     & 5                  & 100                 & ReLU &                      & 3                  & 100                 & \multicolumn{1}{l}{Tanh} &  & 8                  & 50        & 90               \\
		\hline
	\end{tabular}
	\label{table:RKN_ex3}
\end{table}

\begin{table}[!h]
	\centering
	\caption{Hyperparameters of the network structure for the CNF-DPF method.}
	\begin{tabular}{cccclccclccc}
		\hline
		\multicolumn{1}{l}{} & \multicolumn{3}{c}{Encoder}         &                      & \multicolumn{3}{c}{Measurement model}          &  & \multicolumn{3}{c}{NF \& CNF}                                \\ \cline{2-4} \cline{6-8} \cline{10-12} 
		& $N_{\text{layer}}$ & $N_{\text{units}}$ & Act. &                      & $N_{\text{layer}}$ & $N_{\text{units}}$ & Act. &  & $N_{\text{block}}$ & $N_{\text{layer}}$ & $N_{\text{units}}$ \\ \hline
		\cref{subsection:sinusoidal-example}          & 4                  & 64                 & ReLU & \multicolumn{1}{c}{} & 3                  & 32                 & ReLU &  & 2                  & 3                  & 64                 \\
		\cref{subsection:lorenz96_example}          & 4 & 64 & ReLU         & \multicolumn{1}{c}{} & 3                  & 32                 & ReLU &  & 2                  & 3                  & 64                 \\ \hline
	\end{tabular}
	\label{table:DPF_ex1_ex2}
\end{table}

\begin{table}[!h]
    \centering
    \caption{Hyperparameter settings of model training across different methods.}
    \begin{tabular}{cccclccclcc}
    \hline
    \multicolumn{1}{l}{}   & \multicolumn{3}{c}{FBF}                                      &                      & \multicolumn{3}{c}{RKN}                                                                                                              &  & \multicolumn{2}{c}{CNF-DPF}             \\ \cline{2-4} \cline{6-8} \cline{10-11} 
                           & Sec. \ref{subsection:sinusoidal-example}           & Sec. \ref{subsection:lorenz96_example}           & Sec. \ref{subsection:AD_example}           &                      & Sec. \ref{subsection:sinusoidal-example}                                              & Sec. \ref{subsection:lorenz96_example}                                               & Sec. \ref{subsection:AD_example}            &  & Sec. \ref{subsection:sinusoidal-example}           & Sec. \ref{subsection:lorenz96_example}           \\ \hline
    $E_{\text{epoch}}$     & 500               & 5000               & 500               & \multicolumn{1}{c}{} & 500                                                   & 5000                                                   & 500               &  & 500               & 5000               \\
    $N_{\text{batch}}$     & 64                 & 64                & 64                 & \multicolumn{1}{c}{} & 64                                                     & 64                                                    & 64                 &  & 64                 & 64                 \\
    $\mathrm{lr}$ & $5\mathrm{e}^{-4}$ & $5\mathrm{e}^{-4}$ & $5\mathrm{e}^{-4}$ &                      & $1\mathrm{e}^{-3}$ & $5\mathrm{e}^{-4}$ & $5\mathrm{e}^{-4}$ &  & $1\mathrm{e}^{-4}$ & $1\mathrm{e}^{-4}$ \\ \hline
    \end{tabular}
    \label{table:training}
\end{table}

\end{document}